\documentclass{amsart}

\usepackage[hyphens]{url}
\usepackage{hyperref}
\usepackage[hyphenbreaks]{breakurl}
\usepackage{amsmath}
\usepackage{amssymb}
\usepackage{mathrsfs}
\usepackage{url}
\usepackage{tikz}
\usepackage{tikz-cd}
\usepackage[square,numbers]{natbib}
\theoremstyle{plain}
\newtheorem{thm}{Theorem}[section]
\newtheorem*{thm*}{Theorem}
\newtheorem{lem}[thm]{Lemma}
\newtheorem{prop}[thm]{Proposition}
\newtheorem{cor}[thm]{Corollary}
\newtheorem{quest}[thm]{Question}
\newtheorem*{quest*}{Question}
\newtheorem{hyp}[thm]{Hypothesis}

\theoremstyle{definition}
\newtheorem{defn}[thm]{Definition}
\newtheorem{eg}[thm]{Example}

\newtheorem{fact}[thm]{Fact}

\theoremstyle{remark}
\newtheorem{rem}[thm]{Remark}

\def\ind{\mathrel{\raise0.2ex\hbox{\ooalign{\hidewidth$\vert$\hidewidth\cr\raise-0.9ex\hbox{$\smile$}}}}}

\providecommand{\:}{}
\renewcommand{\:}[0]{\ : \ }

\newcommand{\I}[0]{\mathcal{I}}
\newcommand{\U}[0]{\mathcal{U}}

\providecommand{\L}{}
\renewcommand{\L}{\mathcal{L}}

\newcommand{\fhi}[0]{\varphi}
\newcommand{\tleq}{\trianglelefteq}

\newcommand{\res}{{\upharpoonright}}

\renewcommand{\H}{\mathcal{H}}
\newcommand{\C}{\mathcal{C}}
\renewcommand{\I}{\mathcal{I}}

\renewcommand{\P}{\mathcal{P}}
\newcommand{\K}{\mathcal{K}}
\newcommand{\F}{\mathcal{F}}
\renewcommand{\and}{\ \& \ }

\newcommand{\OP}{\mathrm{OP}}
\newcommand{\IP}{\mathrm{IP}}
\newcommand{\NIP}{\mathrm{NIP}}
\newcommand{\SOP}{\mathrm{SOP}}
\newcommand{\NSOP}{\mathrm{NSOP}}
\newcommand{\TP}{\mathrm{TP}}
\newcommand{\NTP}{\mathrm{NTP}}

\newcommand{\PM}{\mathrm{PM}}
\newcommand{\NPM}{\mathrm{NPM}}

\newcommand{\Int}{\mathrm{Int}}
\newcommand{\acl}{\mathrm{acl}}
\newcommand{\Aut}{\mathrm{Aut}}
\title{Toward a new dichotomy in NSOP theories}
\author{Michele Bailetti}
\date{\today}

\begin{document}

\begin{abstract}
 We present a conjectural dichotomy (every $\NSOP$ theory is either $\NSOP_3$ or $\PM^{(2)}$) which can be viewed as a variant of the characterization of stable theories by $\NSOP \cap \NIP$. We provide some partial results and some possible strategies to settle it. If the proposed dichotomy is true, the collapse of the $\NSOP_n$ hierarchy inside $\NTP_2$ follows.
\end{abstract}
\maketitle
\setcounter{tocdepth}{1}
\tableofcontents

%Introduction
\section{Introduction}
\label{introduction}

In classification theory, the tame side of many dividing lines are defined by the absence of local combinatorial properties. In \cite{BAILETTI_2025}, we first gave a general framework to study these kind of dividing lines, using and elaborating the concept of pattern of consistency and inconsistency (which appears in different forms in \cite{shelah1999iunderstandandsay} and \cite{kruckman2023newkim}); then, taking this idea to its limit, we defined and explored several notions of maximal complexity for first-order theories. In particular, we defined a strict hierarchy of properties $(\PM^{(k)})_{1<k<\omega}$ by weakening the notion of \emph{Positive Maximality} ($\PM$) restricting the size of the inconsistency conditions (See \cite{BAILETTI_2025}, Definition 5.2, Definition 6.1).

The weakest property in this hierarchy, $\PM^{(2)}$, is already strong enough to imply (even at the level of formulas) many other combinatorial configurations: it is easy to see that $\PM^{(2)}$ implies $2$-$\TP_2$ (and thus $\IP$, $\OP$, $\TP$) since this property can be described by positive patterns (see \cite{BAILETTI_2025}, Definition 5.1) in which the inconsistency conditions have size $2$. Less trivially (at first look), $\PM^{(2)}$ implies $\SOP_3$ (and thus $\SOP_1/\SOP_2/\TP_1$); indeed, in \cite{kaplan2023generic}, Kaplan, Ramsey and Simon proved a characterization of $\SOP_3$ which can be described by positive patterns (again with inconsistency condition of size $2$). In more details, $\PM^{(2)}$ can be characterized by the existence of a formula which realizes (in some sense) every finite graph while $\SOP_3$, using the aforementioned characterization, is equivalent to the existence of a formula realizing every finite \emph{interval} graph. The tension of realizing every interval graph ($\SOP_3$) while omitting some (non-interval) graph ($\NPM^{(2)}$) seems to suggest the presence of an order ($\SOP$). This observation led to the following 
\begin{quest}\label{question}
Does the following equality hold?
\[\NSOP \cap \NPM^{(2)} \stackrel{?}{=} \NSOP_3\]
In other words, in the context of $\NSOP$ theories, do we have the dichotomy: $\NSOP_3$ or $\PM^{(2)}$?
\end{quest}
Clearly, since both $\SOP$ and $\PM^{(2)}$ individually imply $\SOP_3$, we know that $\NSOP_3 \subseteq \NSOP \cap \NPM^{(2)}$. On the other hand, it's not clear whether there exists an $\NPM^{(2)}$, $\NSOP$ theory with $\SOP_3$ or not. A positive answer to Question \ref{question} would imply the collapse of the $\SOP_n$ hierarchy inside $\NTP_2$, \emph{a central open problem in classification theory} (as mentioned in \cite{mutchnik2026applicationsrealstrictorder}), first asked by Chernikov in \cite{CHERNIKOV-theories-without-TP2}, Problem 3.16 (for $n\geq 3$). Indeed, assuming a positive answer to Question \ref{question}, if there was an $\NSOP_{n+1} \cap \NTP_2$ theory with $\SOP_n$, then it would have $\SOP_3$ (since $\SOP_n$ implies $\SOP_3$) and would be $\NSOP$ (since $\NSOP_{n+1}\subset \NSOP$); this would then imply $\PM^{(2)}$ (by our assumption on the dichotomy question) which, in turn, implies $\TP_2$. Therefore, no $\NTP_2$ theory would then be able to distinguish the $\NSOP_n$ hierarchy. Moreover, in light of recent results by Chernikov (see \cite{chernikov2026sop2equalssop3}) showing that $\NSOP_3 = \NTP_1$, and using Shelah's result establishing $\NTP_1 \cap \NTP_2 = simple$ (see \cite{shelah1990classification}), we get that a positive answer to Question \ref{question} would imply $\NSOP \cap \NTP_2 = simple$. 

In this paper, we focus our attention on Question \ref{question}, analyzing the class of $\NSOP$ theories with $\SOP_3$, following both an inductive approach (concrete examples) and a deductive one (abstract reasoning).

In Section \ref{SOP3-interval-graphs}, we first give a rereading of the proof of the characterization of $\SOP_3$ found in \cite{kaplan2023generic} in terms of realization of interval graphs. Then we show how the proposed dichotomy can be viewed as analogous to the classical dichotomy \emph{Stable} or $\IP$ for $\NSOP$ theories, passing from \emph{coding graphs} to \emph{realizing graphs}. 

In Section \ref{SOP3-NSOP}, we consider various examples of $\SOP_3$ theories omitting $\SOP$ and show the presence of a $\PM^{(2)}$ formula in those cases. 

In Section \ref{Dichotomy-for-free-amalgamation-theories}, we focus our attention on \emph{free amalgamation theories}, defined by Conant in \cite{free-amalgamation_Conant}, generalizing free amalgamation in relational structure. With the small hypothesis of \emph{normality} for the abstract free amalgamation relation, we show that the dichotomy holds for this class of theories. 

Finally, in Section \ref{other-partial-results}, assuming the \emph{real-valued $\NSOP_r$ hierarchy is non-distinct from the integer-valued $\NSOP_n$ hierarchy on sufficiently general grounds} (see \cite{mutchnik2026applicationsrealstrictorder}), we show that we can approximate the dichotomy \emph{as close as possible}; more precisely, we prove that, with the above hypothesis, $\NSOP_n \cap \NPM^{(2)} = \NSOP_3$ for all $2<n<\omega$. Furthermore, assuming the existence of a formula realizing every interval graph ($\SOP_3$) but omitting some specific non-interval graphs (whose absence define the class of interval graph itself), we show how to get a $\SOP$ formula, providing more evidence and yet other strategies to settle the dichotomy question.

%Preliminaries
\section{Preliminaries}
\label{preliminaries}

In what follows, $T$ is a complete first-order theory, $\U \models T$ is a monster model and every subset of it is assumed to be a small. Whenever we declare the existence of some tuples, we mean tuples from the monster model; variables $x,y,z,\cdots$ usually denote tuples of variables, not necessarily singletons. For any $B \subseteq \U$ and (tuple of) variable(s) $x$, $B^x$ denotes the set of tuples from $B$ that can be plugged in for $x$. We often consider \emph{partitioned formulas} $\fhi(x;y)$ in which the free variables are divided into \emph{object} and \emph{parameters} variables ($x$ and $y$ respectively).  
%\subsection{Dividing lines} 
\begin{defn}[see \cite{shelah1990classification,SHELAH-toward-class-unstable-t}] \label{div-lines-def}
    Fix $1<k<\omega$ and $2<n<\omega$. A formula $\fhi(x;y)$ has 
    \begin{enumerate}
        \item the \emph{order property} $(\OP)$ if there are tuples $(a_ib_i)_{i<\omega}$ such that 
        \[\models \fhi(a_i;b_j) \iff i<j\]
        \item the \emph{independence property} $(\IP)$ if there are tuples $(a_i)_{i<\omega}$ and $(b_I)_{I\subseteq \omega}$ such that 
        \[\models \fhi(a_i;b_I) \iff i \in I\]
        \item the \emph{strict-order property} $(\SOP)$ if there are $(b_i)_{i<\omega}$ such that, for all $i<j<\omega$
        \[\fhi(\U;b_i) \subsetneq \fhi(\U;b_j)\]
        \item the \emph{$k$-tree property of the second kind} $(k-\TP_2)$ if there are tuples $(b_{i,j})_{i,j<\omega}$ such that 
        \begin{itemize}
            \item paths are consistent: for every $f \in \omega^\omega$, $\{\fhi(x;b_{i,f(i)}) \: i <\omega\}$ is consistent, and
            \item rows are $k$-inconsistent: for each $i<\omega$, $\{\fhi(x;b_{i,j} \: j<\omega\}$ is $k$-inconsistent (i.e.\ any subset of size $k$ is inconsistent). 
        \end{itemize}
        \item the \emph{$1$-strong order property} $(\SOP_1$) if there are tuples $(b_\eta)_{\eta \in 2^{<\omega}}$ such that:
        \begin{itemize}
            \item for all $f \in 2^\omega$, $\{\fhi(x;b_{f\res i}) \: i <\omega\}$ is consistent
            \item for all $\nu,\eta \in 2^{<\omega}$, if $\nu^\frown\langle0\rangle \tleq \eta$, then $\{\fhi(x;b_{\nu^\frown\langle1\rangle}), \fhi(x;b_\eta)\}$ is inconsistent. 
        \end{itemize}
        \item the \emph{$n$-strong order property} $(\SOP_n)$, for $n\geq 3$, if $|x| = |y|$ and there are tuples $(a_i)_{i<\omega}$ such that 
        \begin{itemize}
            \item for all $i<j<\omega$, $\models \fhi(a_i;a_j)$, and
            \item $\bigwedge_{i<n-1}\fhi(x_i;x_{i+1}) \land \fhi(x_{n-1};x_0)$ is inconsistent. 
        \end{itemize}
    \end{enumerate}
\end{defn}

\begin{defn}[see \cite{BAILETTI_2025}, Definitions 3.1, 3.2, 3.5, 5.1]
    Fix $n<\omega$. An \emph{$n$-pattern} is a pair $(\C,\I)$ with $\C,\I \subseteq \P(n)\times\P(n) \setminus \{(\emptyset,\emptyset)\}$. 

    A formula $\fhi(x;y)$ \emph{realize an $n$-pattern $(\C,\I)$} if there are $(b_i)_{i<n}$ such that
    \begin{itemize}
        \item for every $(A^+,A^-) \in \C$ the following is consistent
        \[\{\fhi(x;b_i) \: i \in A^+\}\cup \{\neg \fhi(x;b_j) \: j\in A^-\}\]
        \item for every $(Z^+,Z^-) \in \I$, the following is inconsistent
        \[\{\fhi(x;b_i) \: i \in Z^+ \}\cup\{\neg \fhi(x;b_j) \: j \in Z^-\}\]
    \end{itemize}
    An $n$-pattern $(\C,\I)$ is 
    \begin{enumerate}
        \item \emph{reasonable} if          \begin{itemize}
            \item[(i)] For all $(Z^+,Z^-) \in \I$ and for all $(Y^+,Y^-) \in \C$, we have that $Z^\epsilon \not \subseteq Y^\epsilon$ for some $\epsilon \in \{+,-\}$. 
            \item[(ii)] For all $(A^+,A^-) \in \C$, $A^+ \cap A^- = \emptyset$.
            \item[(iii)] For all $(Z^+,Z^-) \in \I$, $Z^+ \cap Z^- = \emptyset$. 
        \end{itemize}
        \item \emph{positive} if for every $(X^+,X^-) \in \C \cup \I$, we have that $X^- = \emptyset$. 
    \end{enumerate}
\end{defn}

\begin{defn}[see \cite{BAILETTI_2025}, 6.1] 
    Fix $1<k<\omega$. A formula $\fhi(x;y)$ is $\PM^{(k)}$ if for every $n<\omega$ and every reasonable positive pattern $(\C,\I)$ such that, for every $(Z^+,\emptyset) \in \I$, $|Z^+| = k$, the formula $\fhi$ realizes $(\C,\I)$. 
\end{defn}

\begin{defn}
\label{realization-hypergraphs}
    Fix $1<k<\omega$. A formula $\fhi(x;y)$ \emph{realizes a $k$-hypergraph $H =(V,E)$} if there are tuples $(b_v)_{v \in V}$ such that 
    \begin{itemize}
        \item for every clique $v_0,...,v_{m-1}$, the following is consistent
        \[\{\fhi(x;b_{v_i}) \: i<m\}\]
        \item given pairwise distinct vertices $v_0,...,v_{k-1}$ such that $\neg E(v_0,...,v_{k-1})$, the following is inconsistent 
        \[\{\fhi(x;b_{v_i}) \: i < k\}\]
    \end{itemize}
    In particular, for $k=2$, a formula \emph{realizes a graph} if we can find parameters, indexed by the vertices of the graph, such that, instantiating the formula along cliques gives a consistent set of formulas, while instantiating it on any two vertices which are not connected by an edge gives an inconsistent set of formulas. 
\end{defn}

\begin{prop}[see \cite{BAILETTI_2025}, Proposition 6.3] \label{PMk-char-by-hypergraphs}
    A formula is $\PM^{(k)}$ if and only if it realizes every finite $k$-hypergraph. 

    In particular, a formula is $\PM^{(2)}$ if and only if it realizes every finite graph. 
\end{prop}

%SOP3-interval-graphs
\section{$\SOP_3$ and interval graphs}
\label{SOP3-interval-graphs}

\begin{defn}
    A graph $G = (V,E)$ is called an \emph{interval graph} if for every $v \in V$ there exists a closed interval $I_v$ on $(\mathbb{R},<)$ such that, for every $v\neq w$ in $V$
    \[vEw \iff I_v \cap I_w \neq \emptyset\]
    The set of intervals $\{I_v \: v \in V\}$ is called an \emph{interval representation for $G$ on $(\mathbb{R},<)$}.
    Given an arbitrary linear order $(L,<)$, if the intervals for the representation of $G$ can be chosen from the set of intervals of $(L,<)$, we say that $G$ has an \emph{interval representation on $(L,<)$}.
\end{defn}

\begin{prop}\label{finite-interval-has-rep-in(n,<)}
    For every finite interval graph $G$ there exists $n<\omega$ such that $G$ has an interval representation on $(n,<)$.   
\end{prop}

\begin{proof}
    Let $G$ be a finite interval graph and let $\{I_v \: v \in V\}$ an interval representation of $G$ on $(\mathbb{R},<)$. For each $v \in V$, say $I_v = [l_v,r_v]$ with $l_v \leq r_v$ in \rotatebox{90}{$\mathbb{R}$ }$\ =\mathbb{R}\cup\{-\infty,+\infty\}$.
    The set $E = \{l_v \: v \in V\}\cup\{r_v \: v \in V\} \subset \ $\rotatebox{90}{$\mathbb{R}$ } is finite and has an induced order; thus, since finite linear orders are characterized by size, there exists $n<\omega$ and a order-preserving bijection $f:E \to n$. We claim that $\{[f(l_v),f(r_v)] \: v \in V\}$ is an interval representation for $G$ on $(n,<)$. 

    Let $v\neq w$ in $V$. We have
    \[v E w \iff I_v \cap I_w \neq \emptyset \iff (l_v \leq l_w\leq r_v) \lor (l_w \leq l_v \leq r_w)\]
    \[\iff (f(l_v) \leq f(l_w)\leq f(r_v)) \lor (f(l_w) \leq f(l_v) \leq f(r_w))\]
    \[\iff [f(l_v),f(r_v)] \cap [f(l_w),f(r_w)] \neq \emptyset\]
\end{proof}

\begin{prop}
    $T$ has $\SOP_3$ if and only if there is a formula $\psi(x;y)$ such that, for every $n< \omega$ there are $(b_I)_{I \in \Int(n)}$ such that, for every non-empty $\I \subset \Int(n)$, we have
    \[\bigcap \I \neq \emptyset \iff \models \exists x \bigwedge_{I \in \I} \psi(x;b_I)\]
    where $\Int(n) = \{\ [l,r] \: l \leq r <n\}$ is the set of non-empty intervals on $(n,<)$.
\end{prop}

\begin{proof}
    This proof is similar to the one done in \cite{kaplan2023generic} by Kaplan, Ramsey and Simon. Their proof uses a characterization of $\SOP_3$ given by two formulas; here we review that proof but in terms of the usual definition of $\SOP_3$ (Definition \ref{div-lines-def}, (6) for $n=3$).
    
    Suppose $\fhi(x;y)$ has $\SOP_3$ witnessed by $(a_i)_{i<\omega}$. Fix $n<\omega$ and for $I = [l,r] \in \Int(n) $ define $b_I = a_{3l}a_{3r+2}$.
    We claim that the formula $\psi(x;yz) := \fhi(y,x) \land \fhi(x,z)$ satisfies our condition with $(b_I)_{I \in \Int(n)}$. 

    Let $\I \subseteq Int(n)$. First, assume that $\bigcap \I = \emptyset$; this implies that there are $I_j = [l_j,r_j]$, for $j \in 2$, such that $I_0 \cap I_1 = \emptyset$. Without loss of generality, we can assume $r_0 < l_1$. Assume, towards contradiction that there exists $a$ such that 
    \[a \models \psi(x;b_{I_0})\land \psi(x;b_{I_1})\]
    Then, by definition of $\psi$, we get
    \[\models \fhi(a_{3l_0},a)\land \fhi(a,a_{3r_{0}+2}) \land \fhi(a_{3l_1},a)\land \fhi(a,a_{3r_{1}+2})\]
    Moreover, since $r_0 < l_1$ (in  $\omega$), we have that $3r_0+2 < 3l_1$ and thus, by $\SOP_3$ we have $\models \fhi(a_{3r_0+2},a_{3l_1})$. In particular we have that 
    \[a\ a_{3r_0+2}\ a_{3l_1} \models \fhi(x_0,x_1)\land \fhi(x_1,x_2) \land \fhi(x_2,x_0) \]
    which is supposed to be inconsistent by $\SOP_3$. Thus $\psi(x;b_{I_0}) \land \psi(x;b_{I_1})$ is inconsistent, and, a fortiori, the following is inconsistent
    \[\bigwedge_{I \in \I}\psi(x;b_I)\]

    Assume now that $\bigcap \I \neq \emptyset$. Say $\I = \{I_j \: j<m\} = \{ \ [l_j,r_j] \: j<m \}$.
    Since $\bigcap \I \neq \emptyset$, there exists $k<n$ such that $l_j \leq k \leq r_j$ for all $j<m$. Let $a = a_{3k+1}$. Since $l_j \leq k \leq r_j$, we have that $3l_j < 3k+1 < 3r_j+2$ and thus, by $\SOP_3$ 
    \[\models \fhi(a_{3l_j},a_{k+1}) \land \fhi(a_{k+1},a_{3r_j+2}) = \psi(a_{k+1},b_{I_j})\]
    for every $j<m$. Thus 
    \[\models \exists x\bigwedge_{j<m} \psi(x;b_{I_j})\]

Suppose now that $\psi(x;y)$ satisfies the above condition. Consider the formula
\[\fhi(x_0,x_1,x_2;y_0,y_1,y_2) = \psi(y_0;x_1) \land \psi(x_0;y_2) \land \neg \exists w (\psi(w;y_1)\land\psi(w;x_2))\]
Then $\fhi(\bar{x};\bar{y})\land \fhi(\bar{y},\bar{z}) \land \fhi(\bar{z},\bar{x})$ contains, as a conjunct, the following \[\psi(y_0;x_1)\land \psi(y_0;z_2) \land \neg \exists w (\psi(w;x_1)\land \psi(w;z_2))\]
which is clearly inconsistent. 
It remains to find $(\bar{a}_i)_{i<\omega} = (a_i^0a_i^1a_i^2)_{i< \omega}$ such that, for all $i<j<\omega$, we have $\models \fhi(\bar{a}_i;\bar{a}_j)$. By compactness, it's enough to find, for each $n<\omega$, tuples $(\bar{a}_i)_{i<n}$ with the same property. Fix $n<\omega$; by assumption of $\psi$, find $(b_I)_{I \in \Int(n)}$ and, for each $i \in n$, find $a_i^0$ such that
\[a_i^0 \models \bigwedge \{\psi(x;b_I) \: I \in \Int(n),\  I \ni i\}\]
and let $a_i^1 = b_{[i,n-1]}$ and $a_i^2 = b_{[0,i]}$. Fix $i<j <n$. Now, since $i<j<n$,  $i \in [0,j]$ and thus $\psi(a_i^0;a^2_j)$. Moreover, since $i<j<n$, we have $j \in [i,n-1]$ and so $\psi(a_j^0,a_i^1)$. Finally, again by the fact that $i<j$, we have that $[0,i] \cap [j,n-1] = \emptyset$ and so $\neg \exists z \psi(z;a_j^1)\land\psi(z;a_i^2)$. 
\end{proof}

\begin{rem}\label{rem-SOP3=real-Intervals}
    The above proposition, together with the previous result, shows that $\SOP_3$ is equivalent to the existence of a formula realizing (in the sense of Definition \ref{realization-hypergraphs}) every interval graph. In particular $\PM^{(2)}$ implies $\SOP_3$; one can see this in two ways: first of all, the conditions on $\psi$ in the previous proposition can be described with positive patterns with inconsistency requests of size $2$, and thus by patterns that are realized by a $\PM^{(2)}$ formula. Alternatively, $\PM^{(2)}$ can be characterized by the existence of a formula realizing every finite graph; in particular it realizes every finite interval graph, and thus the theory has $\SOP_3$.  
\end{rem}

We now show that a positive answer to Question \ref{question} would be an analog of the classical dichotomy (either stable or $\IP$, see \cite{shelah1990classification}) in the context of $\NSOP$ theories, passing from \emph{coding graphs} to \emph{realizing graphs}. 

Recall the following definition 
\begin{defn}[\cite{LASKOWSKI-Shelah-codinggraphs}]
A formula $\fhi(y_0,y_1)$ \emph{codes a graph $G=(V,E)$} if there are $\{a_v \: v \in V\}$ such that 
\[\models \fhi(a_v,a_w) \iff vEw\]
A formula \emph{codes graphs} if it codes every graph. 
\end{defn}

\begin{fact}[\cite{LASKOWSKI-Shelah-codinggraphs}, Lemma 2.2]
A theory has $\IP$ if and only if there is a formula that codes graphs. 
\end{fact}

A natural (and thematic in this setting) weakening of the above property is \emph{coding interval graphs}. This turns out to be equivalent to the order property: 

\begin{prop}
    A theory is stable if and only if no formula codes interval graphs.
\end{prop}

\begin{proof}
    Suppose $\fhi(x;y)$ codes interval graphs. For $n<\omega$, consider the graph on $\{v_i \: i<n\}\cup\{w_i \: i<n\}$ where $(w_i)_{i<n}$ is an independent set, $(v_i)_{i<n}$ forms a clique and 
    \[v_i E w_j \iff i<j\]
    It's easy to see that this graph is interval and so, for every $n<\omega$, we can find $(a_ib_i)_{i<n}$ such that 
    \[\models \fhi(a_i,b_j) \iff i<j\]
    Thus $\fhi$ has $\OP$ and the theory is unstable.

    Conversely, if $\fhi(x;y)$ has $\OP$, witnessed by $(a_ib_i)_{i<\omega}$, then $\psi(xy;zw) = \fhi(x,w) $ defines a strict linear order on $(d_i)_{i<\omega}$, where $d_i = a_ib_i$. 
    Now that we have a formula $\psi$ defining a linear order on $(d_i)_{i<\omega}$, consider the formula 
    \[\theta(xy,x'y') = ((\psi(x,x')\lor x=x') \land (\psi(x',y) \lor x' = y))\lor \]\[((\psi(x',x) \lor x = x')\land (\psi(x,y') \lor  x = y'))\]
    This should be intended as follows:
    if $x,y$ and $x',y'$ are the left and right endpoints of two intervals $I = [x,y]$, $I'=[x',y']$, then $I \cap I' \neq \emptyset$ if and only if 
    \[(x\leq x' \leq y) \lor (x' \leq x \leq y')\]
    So $\theta$ describes the condition of non-empty intersection of two intervals. 
    Any finite interval graph has a representation on $(n,<)$ for some $n<\omega$ (see Proposition \ref{finite-interval-has-rep-in(n,<)}); thus, given a finite interval graph $G$, for appropriate $n$, $G$ is an induced subgraph of $\Int(n)$. For $I = [l,r]\in \Int(n)$, let $b_I = d_ld_r$. Then $G$ is coded by $\theta$. 
\end{proof}

\begin{rem}
The above proposition shows that the difference between instability and $\SOP_3$ is simply in the difference between coding and realizing interval graphs, and thus $\NSOP_3$ is a natural extension of the class of stable theories. 
\end{rem}

Question \ref{question} can be answered positively if, assuming the existence of a formula \emph{realizing} every interval graph ($\SOP_3$) but no formula \emph{realizes} every graph ($\NPM^{(2)})$, we can get a formula with $\SOP$.

Exchanging \emph{realizing} with \emph{coding} gives a true statement: 
Assuming the existence of a formula \emph{coding} every interval graph ($\OP$) but no formula \emph{codes} every graph ($\NIP$), then there is a formula with $\SOP$.
Indeed, by the classical dichotomy, an $\NSOP$ theory is either stable or it has $\IP$ (see \cite{shelah1990classification}).

%SOP3-NSOP-theories
\section{Some of those $\SOP_3$ theories omitting $\SOP$}
\label{SOP3-NSOP}

In this section we prove that many examples of $\NSOP$ theories with $\SOP_3$ are $\PM^{(2)}$. For some of the following examples, this fact can be obtained as a consequence of Theorem \ref{freeamwithSOP3isPM^2}, stated and proven in the next section; we include these proofs to provide explicit examples of $\PM^{(2)}$ formulas.

%SOP3 OMITTING SOP4
\subsection{$\SOP_3$ theories omitting $\SOP_4$.}

%GENERIC K_n-FREE GRAPH
\begin{eg}[Generic $K_n$-free graph]
Fix $3\leq n < \omega$ and let $\H = (H,E)$ be the generic $K_n$-free graph (i.e.\ the Fra\"issé limit of the class of finite $K_n$-free graphs; see \cite{Henson-family-cntb-hom-graphs}). Consider the formula 
\[\fhi(\bar{x};y_0y_1) = \bar{x}\cong K_{n-2} \land \bigwedge_{\substack{i<n-2 \\ j \in  2}} x_i E y_j \]

Fix a finite graph $G = (V,e)$ and define $G' = (V',e')$ as follows: 
\[V' = \{v'_i \: v \in V, \ i \in 2\}\]
\[e'(v'_i,w_j') \iff i \neq j \land \neg e(v,w) \land v\neq w\]
Note that $G'$ is bipartite; in particular $G'$ is $K_n$-free and, by universality, we can embed it in $\H$. We claim that $(v'_0v'_1)_{v \in V}$ (viewed as embedded in $\H$) witness the realization of $G$ in $\H$ by $\fhi$.

First, let $v\neq w \in V$ such that $\neg  e(v,w)$. Then $v'_0 E w'_1$. If there was $\bar{a} \cong K_{n-2}$ such that $a_i E v'_j \land a_i E w'_k$ for $i<n-2$ and $j,k \in 2$, we would have that $\bar{a}v'_0w'_1 \cong K_n$. Contradiction. Thus $\fhi(\bar{x};v'_0v'_1) \land \fhi(\bar{x};w'_0w'_1)$ is inconsistent.   

On the other hand, assume $C$ is a clique in $G$. Then $(v'_0v'_1)_{v \in C}$ is an independent set; therefore, by the extension property of $\H$, we can find $\bar{a}$ such that 
\[\bar{a} \models \bigwedge_{v \in C}\fhi(\bar{x};v'_0v'_1)\]

In conclusion, $\fhi(\bar{x},y_0y_1)$ realizes every finite graph and thus, by Proposition \ref{PMk-char-by-hypergraphs}, for $k=2$, it is $\PM^{(2)}$. 
\end{eg}

%HENSON DIGRAPHS
\begin{eg}[Henson digraphs]
Let $\F$ be a finite set of tournaments (each of size at least $3$) and consider the class $\K_\F$ of all $\F$-free directed graphs; it can be shown that $\K_\F$ is a Fra\"issé class (see
\cite{Henson_cntbl-hom-rel-structures}). Let $T_\F$ be the complete theory of the Fra\"issé limit $\H_\F$ of $\K_\F$.
We claim that $T_\F$ is $\PM^{(2)}$.  

Let $k := \min\{|F| \: F \in \F\}$ and let $F^* \in \F$ such that $|F^*| = k$. Fix $v^*,w^* \in F^*$ and let $B$ be the digraph obtained from $F^*$ by removing the edge (without loss of generality) $R(v^*,w^*)$ between $v^*$ and $w^*$. Consider the formula 
\[\fhi(\bar{x};y_0y_1) = \bar{x}y_0y_1 \cong B\]
where $y_0$ and $y_1$ represent $v^*$ and $w^*$ respectively. 

We claim that $\fhi$ is $\PM^{(2)}$. 
Fix a finite graph $G = (V,E)$ and define the digraph $G' = (V',R')$ as follows:
\[V' = \{v'_i \: v \in V, \ i \in 2\}\]
\[R(v'_i,w'_j) \iff i = 0 \land j = 1 \land v \neq w \land \neg vEw\]

If we view $G'$ as a graph (forgetting the direction of edges), it is bipartite and thus, as a digraph, it contains no tournament. By universality, we can embed $G'$ in $\H_\F$. 
If $v\neq w \in G$ are such that $\neg v Ew$, then, in $G'$ we have $R(v'_0,w'_1)$ and so it in inconsistent to find $\bar{a}$ such that $\bar{a}v'_0v'_1 \cong B$ and $\bar{a}w'_0w'_1 \cong B$ since, otherwise, we would get that $\bar{a}v'_0w'_1 \cong F^*$. 
On the other hand, if $C$ is a clique in $G$, then $(v'_0v'_1)_{v \in C}$ is a set of vertices without any directed edge; thus, by the extension properties of $\H_\F$, we can find $\bar{a}$ witnessing the consistency of 
\[\bigwedge_{v \in C}\fhi(\bar{x};v'_0v'_1)\]
Indeed, the largest tournament in $\bar{a}(v_0'v_1')_{v\in C}$ has size $k-1$ and thus it is not in $\F$. 
In conclusion $\fhi$ is $\PM^{(2)}$. 
\end{eg}

%BOWTIE-FREE GRAPH
\begin{eg}[Universal bowtie-free graph, \cite{patel2006family}]
Let $T$ be the theory of the countable universal existentially closed bowtie-free (i.e. $(K_3+K_3)$-free) graph $\U_3$. Consider the formula 
\[\fhi(x_0x_1x_2;y_0y_1) = `` x_0x_1x_2 \cong K_3" \land x_0 E y_0 \land x_0 E y_1 \land \bigwedge_{i,j} x_i \neq y_j\]

Fix a finite graph $G = (V,E)$. Define $G' = (V',E')$ as follows: 
\[V' = \{v_i' \: v \in V, i \in 2\}\] \[v_i' E' w_j' \iff i \neq j \land \neg v E w\land v\neq w\] Note that $G'$ is triangle-free and thus bowtie-free. By universality, we can find a copy of $G'$ in $\U_3$. We claim that $(v_0'v_1')_{v \in V}$ witnesses the realization of $G$ by $\fhi$. 

If $\neg v E w$ in $G$, then there is an edge between $v_0'$ and $w_1'$; thus there are no $a_0a_1a_2 \cong K_3$ such that $\models \fhi(a_0a_1a_2;v_0'v_1')\land \fhi(a_0a_1a_2;w_0'w_1')$. Otherwise we would get a bowtie.

Moreover, if $C$ is a clique in $G$, we can extend $\U_3$ by a triangle $a_0a_1a_2$ where $\bigwedge_{v \in C}a_0Ev_0' \land a_0v_1'$, and no other edges. Since $(v_0'v_1')_{v \in C}$ is an independent set, this is a bowtie-free graph. Since $\U_3$ is existentially closed in the class of bowtie-free graphs, we can find a realization $a_0'a_1'a_2'$ in $\U_3$. Thus 
\[\bigwedge_{v \in C}\fhi(\bar{x};v'_0v'_1)\]
is consistent.
In conclusion, $\fhi$ is $\PM^{(2)}$.
\end{eg}

%GENERIC c-NILPOTENT GROUP OF EXP P>c>2
\begin{eg}[Generic $c$-nilpotent group of exponent $p>c>2$, \cite{Mod-the-prop_c-nilpotent-groups}]
Let $\L_{c,\mathbb{F}}$ be the language of $\mathbb{F}$-vector spaces expanded by a binary function symbol $[-,-]$ and predicates $(P_i)_{1 \leq i \leq c+1}$.  
We prove $\PM^{(2)}$ for the $\L_{c,\mathbb{F}}$-theories $T_{c,\mathbb{F}}$ ($c\geq 3$) of Fra\"issé limits of Lazard Lie Algebras ($\mathrm{LLA}$) over a finite field and, via Lazard correspondence, deduce $\PM^{(2)}$ for the associated groups. 
It suffices to show that $T_{3,\mathbb{F}}$ is $\PM^{(2)}$ since each witness of the realization of a finite graph $G$ is (will be) in a $3$-nilpotent $\mathrm{LLA}$ which may be regarded as a $c$-nilpotent $\mathrm{LLA}$ for any $c\geq 3$ with the interpretation of the Lazard series $P_i = 0$ for all $c>3$ and thus may be embedded in a model of $T_{c,\mathbb{F}}$. 

Define 
\[\fhi(x;yy',zz') = [x,y]=y' \land [x,z]=z'\]
Given a finite graph $G$, define an $\mathbb{F}$-vector space $V$ with basis \[\{a_v,a_v',b_v,b_v',d_{v,w} \: v\neq w \in G, \neg v Ew\}\]
and, for $vEw$, let $d_{v,w}= 0$. Define $[-,-] \colon V^2 \to V$ by $[a_v',b_w] = d_{v,w} $ and $[x,y] = 0$ for all $x,y$ such that $\{x,y\} \neq \{a_v',b_w\}$. One can show that $(V,[-,-])$ is a $2$-nilpotent Lie algebra. We give $V$ a flag structure by interpreting $P_1,P_2,P_3,P_4$ as 
\[P_1(V) = V  \supseteq P_2(V) = \langle a_v',b_w',d_{v,w} \: v,w \in G\rangle \]\[\supseteq P_3(V) = \langle d_{v,w} \: v,w \in G \rangle \]
and $P_4(V) = 0$; one can check that $[P_i(V),P_j(V)] \subseteq P_{i+j}(V)$. We can thus view $V$ as a substructure of $\U \models T_{3,\mathbb{F}}$.

We claim that $\fhi$ realizes $G$ witnessed by $(a_va_v'b_vb_v')_{v \in G}$.
Fix $v,w \in G$ such that $\neg v E w$ and suppose, toward contradiction, that there is $d \models \fhi(x;a_va_v'b_vb_v') \land \fhi(x;a_wa_w'b_wb_w')$. We have 
\[0 \neq d_{v,w} = [a_v',b_w] = [[d,a_v],b_w]\]
Using the Jacobi identity, we get
\[= [d,[a_v,b_w]] - [a_v,[d,b_w]] = [d,0] - [a_v,b_w'] = - [a_v,b_w'] = 0\]
A contradiction. 
On the other hand, suppose $C$ is a clique in $G$. Let $W = \{W,[-,-],W_1,...,W_4\}$ be the substructure generated by $ \{ a_va_v'b_vb_v' \: v \in C\}$. Thus $W$ is the span of these vectors together with a trivial Lie bracket. 
\begin{lem}[see \cite{Mod-the-prop_c-nilpotent-groups}, Lemma 5.7]
   \label{lemma-nilpotentliealgebra} 
   Let $\mathbb{F}$ be a field. Suppose $V,W$ are $\mathbb{F}$-vector spaces and $[-,-]_0 : V^2 \to W$ is an alternating bilinear map. Define $[-,-]: (V\oplus W)^2 \to (V \oplus W)$ as 
   \[[(v,w),(v',w')] = (0,[v,v']_0)\]
   Then $[-,-]$ gives $V \oplus W$ the structure of a $2$-nilpotent Lie algebra over $\mathbb{F}$. 
\end{lem}
Let $W_*$ be an $\mathbb{F}$-vector space spanned by $\{a_vb_v \: v \in C\} \cup \{c_*\}$, where $c_*$ is a new basis element. Let $W_{**}$ be the $\mathbb{F}$-vector space spanned by $\{a_v'b_v' \: v \in C\}$. Define $[-,-]_* : W_* \to W_{**}$ by
\[[c_*,a_v] = a_v'\]
\[[c_*,b_v] = b_v'\]
\[[x,y] = 0 \text{ for all other $x,y$}\]
By Lemma \ref{lemma-nilpotentliealgebra}, let $\Tilde{W} = W_* \oplus W_{**}$ with $\Tilde{[-,-]}$ the induced Lie algebra. We may view $\Tilde{W}$ as an extension of $W$, with the flag structure defined by $P_1(\Tilde{W}) = \Tilde{W}$ and $P_i(\Tilde{W}) = P_i(W)$ for $i=2,3,4$. Thus we may embed $\tilde{W}$ in $\U \models T_{c,\mathbb{F}}$ and thus get that
\[\bigwedge_{v \in C} \fhi(x;a_va_v'b_vb_v')\]
is consistent as witnessed by $c_*$. 
In conclusion, $\fhi$ is $\PM^{(2)}$. 
\end{eg}

%SOPN OMITTING SOPN+1
\subsection{$\SOP_n$ theories omitting $\SOP_{n+1}$}

%UNIVERSAL GRAPH OMITTING ODD CYCLES OF LENGTH \LEQ n
\begin{eg}[Universal graph omitting odd cycles of length $\leq n$]
Let $T$ be the complete theory of $\U_n$, the universal and existentially closed countable graph omitting odd cycles of length $\leq n$ (see \cite{Universal-graph-no-odd-cycles-leqn,SHELAH-toward-class-unstable-t}). 
We can assume $n$ is odd (otherwise work with $n-1$). 
Consider the formula
\[\fhi(\bar{x};yz) = \bar{x} \cong P_{n-2} \land e(y,x_0)\land e(x_{n-3},z) \land (\text{no other edges})\]
where $e$ is the edge relation and $P_{n-2}$ is a path with $n-2$ vertices. 
We claim that $\fhi$ is $\PM^{(2)}$.

Fix a finite graph $G = (V,E)$ and define a graph $G' = (V',E')$ as follows: 
\[V' = \{v'_i \: v \in V, i <2\} \]
and 
\[v_i' E' w'_j \iff i\neq j \land \neg v E w \land v\neq w\]

Being bipartite, $G'$ has no odd cycle and thus it embeds by universality in $\U_n$ (maintaining the same names for the vertices). 

Consider $v\neq w$ in $G$ with $\neg v E w$. If $\fhi(\bar{x},v'_0v'_1) \land \fhi(\bar{x},w'_0w'_1)$ was realized by some $\bar{a}$, we would have an $n$-cycle $v'_0 a_0 ... a_{n-3}w'_1$. 

Given a clique $C$ in $G$, extend $\U_n$ with vertices $a_C^0...a_C^{n-3}=\bar{a}_C$ such that $\bar{a}_C \cong P_{n-2}$ and
\[\bigwedge_{v \in C} e(a_C^0, v'_0) \land e(a_C^{n-3},v'_1)\]
and no other edges. If $\U_n \cup \bar{a}_C$ had an odd cycle $O$ of length $\leq n$, then, since $\U_n$ has no odd cycle of length $\leq n$, $O$ must involve some vertices from $\bar{a}_C$; but $|\bar{a}_C| = n-2$ and there's no other edge apart from those to $\{v'_i \: v \in C, i<2\}$ which form an independent set. Thus $O$ must involve $\bar{a}_C,\{v'_i \: v \in C, i<2\}$ and some other vertices connecting two points of the independent set. But then $O$ would have more then $n$ vertices. Therefore $\U_n \cup \bar{a}_C$ has no odd cycle of length $\leq n$. Since in $\U_n \cup \bar{a}_C$ we have that 
\[\bar{a}_C \models  \bigwedge_{v \in C}\fhi(\bar{x};v'_0v'_1)\]
by existential closure of $\U_n$, we can find $\bar{a}_C' \in \U_n$ witnessing the consistency of 
\[\bigwedge_{v \in C}\fhi(\bar{x},v_0'v_1')\]

Thus $\{v'_i \: v \in V, i<2\}$ witness the realization of $G$ by $\fhi$. 
In conclusion, $T$ is $\PM^{(2)}$. 
\end{eg}

%UNIVERSAL DIRECTED GRAPH OMITTING CYCLES OF LENGTH \LEQ n
\begin{eg}[Universal directed graph omitting cycles of length $\leq n$]
Let $T$ be the theory of $\vec{\U}_n$, the universal and existentially closed countable directed graph omitting cycles of length $\leq n$ (see \cite{SHELAH-toward-class-unstable-t}). Consider the fomula 
\[\fhi(\bar{x},yz) = \bar{x} \cong \vec{P}_{n-2} \land e(x_0,y)\land e(z,x_{n-3}) \land (\text{no other directed edges)}\]
where $e$ directed edge relation and $\vec{P}_{n-2}$ is a directed path with $n-2$ vertices. We claim that $\fhi$ is $\PM^{(2)}$. 

Fix a finite graph $G = (V,E)$ and define a directed graph $G' = (V',E')$ as follows:
\[V' = \{v_i' \: v \in V,\ i<2\}\]
and 
\[v_i'E'w_j' \iff i=1 \land j=0 \land \neg v E w \land v\neq w\]
Notice that $G'$ is bipartite and the direction of the edges are always from the set of vertices with index $0$ to those with index $1$. Therefore $G'$ has no directed cycle and thus it embeds into $\vec{\U}_n$ by universality. 

Consider $v\neq w$ in $G$ with $\neg v E w$. If $\fhi(\bar{x},v_0'v_1') \land \fhi(\bar{x},w_0'w_1')$ was realized by some $\bar{a}$, then $v_1'a_0...a_{n-3}w'_0$ would be a directed $n$-cycle. Contradiction. 

On the other hand, let $C$ be a clique in $G$. Extend $\vec{\U}_n$ with vertices $a_C^0...a_C^{n-3}=\bar{a}_C$ such that $\bar{a}_C \cong \vec{P}_{n-2}$ and 
\[\bigwedge_{v \in C}e(v_0',a_C^0) \land e(a_C^{n-3},v_1')\]
and no other edges. If $\vec{\U}_n$ had a directed cycle $O$ of length $\leq n$, then, since $\vec{\U}_n$ has no such cycles, $O$ must involve some vertices from $\bar{a}_C$; but $|\bar{a}_C| = n-2$ and there is no other directed edge apart form those with $\{v_i' \: v \in C\}$ which form an independent set. Thus $O$ must involve $\bar{a}_C$, $\{v_i \: v \in C\}$ and some other vertices connecting two points from the independent set. But then $O$ would have more then $n$ vertices. Therefore $\vec{\U}_n$ has no odd cycle of length $\leq n$. Since in $\vec{\U}_n \cup \bar{a}_C$ we have that \[\bar{a}_C \models \bigwedge_{v \in C}\fhi(\bar{x},v_0'v_1')\]
by existential closure of $\vec{U}_n$, we can find $\bar{a}_C'$ in $\vec{\U}_n$ witnessing the consistency of \[\bigwedge_{v \in C} \fhi(\bar{x},v_0'v_1')\]
Thus $\{v_i' \: v \in V,\ i<2\}$ witness the realization of $G$ by $\fhi$. In conclusion, $T$ is $\PM^{(2)}$. 
\end{eg}

%FREE nTH ROOT OF THE COMPLETE GRAPH
\begin{eg}[Free $n$th root of the complete graph; see \cite{CASANOVAS-WAGNER_Free-roots-Complete-graphs}]

Let $\L = \{R_1,...,R_n\}$ and consider the theory $T_n^-$ of all $n$-colored graphs $R_1,...,R_n$ such that $R_n$ is total and $R_iR_j \subseteq R_{\min\{i+j,n\}}$ for all $0<i,j<n$. Then, let $T_n$ be the model-completion of $T_n^-$. 
As a notation, let $R_0$ be equality. Another way of looking at this theory is that of a metric space of diameter $n$ with distances in $\{0,...,n\}$; the identification is made by 
\[d(x,y) = \min \{i \leq n \: x R_i y\}\]
\begin{fact}[See \cite{CASANOVAS-WAGNER_Free-roots-Complete-graphs}, Proposition 1] \label{factaboutFreenthroot}
    For all finite disjoint $A_1,...,A_{n}$, there exists a point of distance $i$ to all points in $A_i$ for $1 \leq i \leq n$, unless it violates some instance of the triangle inequality. 
\end{fact}
Consider the formula 
\[\fhi(x;y_0y_1) = x R_1 y_0 \land \bigwedge_{i<3}\neg x R_i y_1\]
i.e.\ $d(x,y_0) = 1 \land d(x,y_1) > 2$. We claim that $\fhi$ is $\PM^{(2)}$. 

Given a finite graph $G = (V,E)$, define $\bar{b}_v = b_v^0b_v^1$ such that $d(b_v^0,b_v^1) = 2$ and, for all $v,W \in V$, if $vEw$ then $d(b_v^i,b_w^j) = 2$ for $i,j < 2$ but if $\neg v E w$ then $d(b_v^i,b_w^j) = 1$. 
First, we show that no triangle inequality fails: 
fix distinct $v,w,u \in V$ and $i,j,k<2$. We have 
\[d(b_v^i,b_u^j) \leq 2 = 1 + 1 \leq d(b_v^i,b_w^j) + d(b_w^j,b_u^k)\]
Now, if $\neg v E w$ in $G$, we have that 
\[\fhi(x;\bar{b}_v) \land \fhi(x;\bar{b}_w)\]
is inconsistent. Indeed, if there is $a \models \fhi(x;\bar{b}_v) \land \fhi(x;\bar{b}_w)$, then 
\[2 < d(a,b_w^1) \leq d(a,b_v^0) + d(b_v^0,b_w^1) = 1 + 1 = 2\]
A contradiction. 
On the other hand, if $C$ is a clique in $G$, all the points $(b_v^0b_v^1)_{v \in C}$ are equidistant (of distance $2$) from each others and thus, using Fact \ref{factaboutFreenthroot}, we get
\[\models \exists x \bigwedge_{v \in C} \fhi(x;\bar{b}_{v})\]    
\end{eg}

%SOPn for all n but NSOP
\subsection{An $(\SOP_n)_{n<\omega}$ theory omitting $\SOP$}

%RATIONAL URYSHON SPHERE
\begin{eg}[Rational Uryshon sphere. See, e.g., \cite{CONANT-TERRY_Urysohn-sphere}]
    The Urysohn sphere $\mathbb{U}$ is the unique universal and homogeneous complete separable metric space with diameter $1$. It can be view as the completion of the Fra\"issé limit ($U$, often called the \emph{rational} Urysohn sphere) of the class of finite metric spaces with rational distances and diameter $\leq 1$ in the language $\{d(x,y) \leq r \: r \in \mathbb{Q}\cap [0,1]\}$. We prove that the $U$ is $\PM^{(2)}$. We remark that it is possible to show that $\mathbb{U}$ is an elementary extension of $U$ and so $\mathrm{Th}(\mathbb{U})$ is also $\PM^{(2)}$.

    Consider the formula 
\[\fhi(x;y_0y_1) = d(x,y_0) \leq \frac{1}{3} \land d(x,y_1) > \frac{2}{3}\]
We claim that $\fhi$ is $\PM^{(2)}$. 

Given a finite graph $G = (V,E)$, for all $v \in V$ define $\bar{b}_v = b_v^0b_v^1$ where $d(b_v^0,b_v^1) = \frac{2}{3}$ and, for distinct $v,w \in V$
\[d(b_v^i,b_w^j) = \begin{cases}
    \frac{1}{3} \quad \text{ if } i\neq j \land \neg v E w\\
    \frac{2}{3} \quad \ \text{ otherwise}
\end{cases}\]
For all distinct $v,w,u \in V$ and for all $i,j,k <2$, we have
\[d(b_v^i,b_u^k) \leq \frac{2}{3} = \frac{1}{3}+\frac{1}{3} \leq d(b_v^i,b_w^j)+d(b_w^j,b_u^k)\]

So $\{b_v^i \: v \in V, \ i\in 2\}$ is a finite metric space with diameter $\leq 1$ and we can embed it in $U$ by universality.

Now, suppose $\neg vEw$. If $\fhi(x;\bar{b}_v) \land \fhi(x;\bar{b}_w)$ is consistent (witnessed by $a$), we have 
\[\frac{2}{3} < d(a,b_w^1) \leq d(a,b_v^0)+d(b_v^0,b_w^1) \leq \frac{1}{3} + \frac{1}{3} = \frac{2}{3}\]
A contradiction. Thus $\fhi(x;\bar{b}_v) \land \fhi(x;\bar{b}_w)$ is inconsistent. On the other hand, if $C$ is a clique in $G$, we can take $a$ such that, for each $v \in C$, $d(a,b_v^0) = \frac{1}{3}$ and $d(a,b_v^1) = 1$. It's easy to see that no triangle inequality is violated and so
\[\models \exists x \bigwedge_{v \in C}\fhi(x;\bar{b}_v)\]
    
\end{eg}

%free-amalgamation-theories
\section{Dichotomy for free amalgamation theories}
\label{Dichotomy-for-free-amalgamation-theories}

In this section, we restrict our attention to \emph{free amalgamation theories}, defined in \cite{free-amalgamation_Conant} by Conant with the idea of abstracting free amalgamation in relational structures. We show that, in this restricted framework, the dichotomy holds. 

Recall the following 
\begin{defn}
    A tuple $a \in \U$ is \emph{algebraic over $B$} if it is in some finite set definable with parameters from $B$ (i.e.\ there is a formula $\fhi(x,y)$ and a tuple $b \in B^y$ such that $\fhi(\U,b)$ is finite and $a \models \fhi(x,b)$).  

    Given a set $B$, the \emph{algebraic closure of $B$} is the set $\acl(B)$ of all $a$ that are algebraic over $B$. 

    In what follows, a set $B$ is \emph{closed} if $B = \acl(B)$. 
\end{defn}

\begin{defn}\label{free-amalg-axioms}
    For a ternary relation $\ind$ on subsets of $\U$, define the following axioms: 
    \begin{enumerate}
        \item \emph{Invariance}: For all $A,B,C$, if $A \ind_C B$ and $\sigma \in Aut(\U)$, then $\sigma(A) \ind_{\sigma(C)}\sigma(B)$. 
        \item \emph{Monotonicity}: For all $A,B,C$, if $A \ind_C B$ and $A_0 \subseteq A$, $B_0 \subseteq B$, then $A_0 \ind_C B_0$. 
        \item \emph{Symmetry}: For all $A,B,C$, if $A\ind_C B$ then $B \ind_C A$. 
        \item \emph{Full transitivity}: For all $A$ and $D \subseteq C \subseteq B$, $A \ind_D B$ if and only if $A \ind_C B$ and $A \ind_D C$. 
        \item \emph{Full existence}: For all $B,C$ and tuple $a$, if $C$ is closed then there is $a' \equiv_C a$ such that $a' \ind_C B$. 
        \item \emph{Stationarity}: For all closed $C$ and closed tuples $a,a',b$ with $C \subseteq a \cap b$, if $a \ind_C b$ and $a'\ind_C b$ and $a' \equiv_C a$, then $ab \equiv_Ca'b$.
        \item \emph{Freedom}: For all $A,B,C,D$, if $A \ind_C B$ and $C \cap AB \subseteq D \subseteq C$, then $A \ind_D B$.
        \item \emph{Closure}: For all closed $A,B,C$, if $C \subseteq A \cap B$ and $A \ind_C B$, then $AB$ is closed.  
        
    \end{enumerate}
\end{defn}
\begin{defn}\label{normality-def}
    A ternary relation $\ind$ satisfies \emph{normality} if, for all $A,B,C$, we have that $A \ind_C B$ implies $AC \ind_C B$.
\end{defn}
\begin{rem}
    In \cite{free-amalgamation_Conant}, Conant defines a free amalgamation relation as a ternary relation $\ind$ satisfying all the axioms in Definition \ref{free-amalg-axioms} and a free amalgamation theory as a theory with a free amalgamation relation. The property of normality (Definition \ref{normality-def}) is not included in Conant's definition of free amalgamation theory but, in every example that we know, the relation satisfies normality. Thus, in our definition, we include the normality axiom. 
\end{rem}

\begin{defn}\label{free-ama-def}
A ternary relation on subsets of $\U$ is a \emph{free amalgamation relation} if it satisfies all the axioms in Definition \ref{free-amalg-axioms} and normality (Definition \ref{normality-def}). A theory is a \emph{free amalgamation theory} if it has a free amalgamation relation.    
\end{defn}

In \cite{free-amalgamation_Conant}, generalizing a result of Patel (\cite{patel2006family}), Conant proved the following 
\begin{fact}[\cite{free-amalgamation_Conant}, Theorem 4.4]
    If $T$ is a free amalgamation theory, then $T$ is $\NSOP_4$. 
\end{fact}

So, assuming free amalgamation, we restrict our attention to those $\NSOP$ theories that omit $\SOP_4$ and have a free amalgamation relation.   

\begin{defn}[See, e.g., \cite{free-amalgamation_Conant}, Definition 7.5]
    Let $\ind$ be a ternary relation on subsets of $\U$, $\mu$ an ordinal, $(b_i)_{i < \mu}$ a sequence of tuples from $\U$ and $C \subseteq b_0$. Then, we say that $(b_i)_{i<\mu}$ is \emph{$\ind$-independent over $C$} if for all $i<\mu$
    \begin{enumerate}
        \item $b_i \equiv_C b_0$
        \item $b_i \ind_C b_{<i}$
    \end{enumerate}
\end{defn}

\begin{lem}[\cite{free-amalgamation_Conant}, Lemma 7.6]\label{dividingalongIndepSequences}
    Suppose $\ind$ is a free amalgamation relation for $T$. Fix closed tuples $a,b$ and let $C = a\cap b$. Suppose, for some ordinal $\alpha$, $(b_i)_{i<\alpha}$ is $\ind$-independent over $C$, with $b_0 = b$. Then, there is $a^*$ such that $a^*b_i \equiv_C ab$ for all $i<\mu$. 
\end{lem}

The following lemma is the only place where we need normality.

\begin{lem}\label{pairsofIndeptuples}
    Let $\ind$ a free amalgamation relation. Suppose $(b_i)_{i<\omega}$ is a  $\ind$-independent sequence over $C$ and each $b_i$ is closed. Then $(b_{2i}b_{2i+1})_{i<\omega}$ is $\ind$-independent over $C$. 
\end{lem}
\begin{proof}
First, for $i<\omega$, we show that $b_{2i}b_{2i+1} \equiv_C b_0b_1$. 

Indeed, by assumption, $b_{2i} \equiv_C b_0$; thus, $b_{2i}b_{2i+1}\equiv_C b_0\sigma(b_{2i+1})$ for some $\sigma \in \Aut(\U/C)$. 
Now, $b_1 \ind_C b_0$ and, since $b_{2i+1}\ind_C b_{2i}$, by invariance, we get $\sigma(b_{2i+1}) \ind_C b_0$. Moreover, clearly, $\sigma(b_{2i+1})\equiv_C b_1$. Therefore, by stationarity, we have that $\sigma(b_{2i+1})b_0 \equiv_C b_1b_0$. In conclusion
\[b_{2i}b_{2i+1}\equiv_C b_0\sigma(b_{2i+1}) \equiv_C b_0b_1\]

We now show that, for every $i<\omega$, we have 
\[b_{2i}b_{2i+1} \ind_C (b_{2j}b_{2j+1})_{j<i}\]

By assumption we have 
\[b_{2i+1} \ind_C b_{<2i+1}\]
and thus 
\[b_{2i+1}\ind_C b_{2i}(b_{2j}b_{2j+1})_{j<i}\]
By full transitivity, since $C \subseteq b_{2i} \subseteq b_{2i}(b_{2j}b_{2j+1})_{j<i}$, we get
\[b_{2i+1}\ind_{b_{2i}} b_{2i}(b_{2j}b_{2j+1})_{j<i}\]
By monotonicity, we have
\[b_{2i+1}\ind_{b_{2i}}(b_{2j}b_{2j+1})_{j<i}\]
and by normality we have
\[b_{2_i}b_{2i+1}\ind_{b_{2i}} (b_{2j}b_{2j+1})_{j<i}\]
Finally, since $b_{2i} \ind_C (b_{2j}b_{2j+1)_{j<i}}$, we can conclude by  symmetry and full transitivity
\[b_{2i}b_{2i+1} \ind_C (b_{2j}b_{2j+1})_{j<i}\]
\end{proof}

\begin{lem}[\cite{free-amalgamation_Conant}, Lemma 7.15] \label{free-am-implies-alg-ind} Suppose $\ind$ satisfies invariance, monotonicity, full existence and stationarity. Then, for any closed tuples $a,b$ and closed set $C$, with $C \subseteq a \cap b$, if $a \ind_C b$ then $a \cap b = C$.
    
\end{lem}

We are now ready to prove the main theorem of this section. 

\begin{thm}\label{freeamwithSOP3isPM^2}
    Let $T$ be a free amalgamation theory. Then $T$ has $\SOP_3$ if and only if $T$ is $\PM^{(2)}$. 
\end{thm}

\begin{proof}
    First, if $T$ is $\PM^{(2)}$ then $T$ has $\SOP_3$ (and this holds even without assuming free amalgamation).

    Suppose now that $T$ is a free amalgamation theory with $\SOP_3$. 
     Let $\ind$ be the free amalgamation relation and let $\psi(x;y)$ be a formula witnessing $\SOP_3$ with $(b_i)_{i<\omega}$. So, for all $i<j$ we have $\psi(b_i,b_j)$ and 
    \[\models \neg \exists x_0x_1x_2(\psi(x_0,x_1)\land \psi(x_1,x_2)\land \psi(x_2,x_0))\]
    For every $i<\omega$, let $\Bar{b}_i := \acl(b_i)$. Without loss of generality, we can assume that $(\Bar{b}_i)_{i<\omega}$ is indiscernible over $C := \Bar{b}_0 \cap \Bar{b}_1$ (given an indiscernible witness for $\SOP_3$, consider an indiscernible realization of the EM-type of the sequence of the algebraic closure of each tuple in the sequence and map the part corresponding to the original sequence back by an automorphism; this automorphism will send the algebraic closure of each tuple of the new sequence to the algebraic closure of the tuples of the original sequence, maintaining indiscernibility; finally this indiscernible sequence will be indiscernible over the common intersection).  

    We claim that the formula 
    \[\fhi(x;yz) = \psi(x,y)\land \psi(z,x)\]
    is $\PM^{(2)}$. 

    The idea is the following: 
    To realize every finite graph, each vertex $v$ will be represented by a pair $(b_v^0b_v^1)$ and the formula says 
    \begin{center}
    \begin{tikzpicture}[node distance = 15mm]
        \node[] (b0) {$b_v^0$};
        \node[] (b1)[below of =b0]{$b_v^1$};
        \node[] (x)[left of =b0]{$x$};

        \draw[->] (x) -- (b0);
        \draw[->] (b1)--(x);
    \end{tikzpicture}
    \end{center}
    Using the axioms of a free amalgamation relation and the tuples witnessing $\SOP_3$, we build cliques independently enough to apply Lemma \ref{dividingalongIndepSequences} and get consistency, and non-edges $\neg v E w$ in such a way that we force a three-cycle (inconsistent by $\SOP_3$):
    \begin{center}
    \begin{tikzpicture}[node distance = 15mm]
        \node[] (b0v) {$b_v^0$};
        \node[] (phantom) [right of =b0v]{};
        \node[] (b1v)[below of =b0]{$b_v^1$};
        \node[] (b0w) [right of = phantom]{$b_w^0$};
        \node[] (b1w) [below of = b0w]{$b_w^1$};
        \node[] (x)[above of =phantom]{$x$};

        \draw[->] (b1v) to [out=135](x);
        \draw[->,thick,red] (x) to (b0v);
        \draw[->,thick,red] (b1w) to [out=45,in=45] (x);
        \draw[->] (x) to (b0w);
        \draw[->,thick,red] (b0v) to (b1w);

    \end{tikzpicture}
    \end{center}

    By induction on $n<\omega$, we show that, if $G = (n,E)$ is a finite graph, we can find $(b_i^0b_i^1)_{i<n}$ such that 
    \begin{enumerate}
        \item[(i)] For all $i<n$, $b_i^0 \ind_C b^0_{<i}b^1_{<i}$.
        \item[(ii)] For all $i<j< n$ such that $\neg i Ej$, we have $b_i^0b_j^1 \equiv_C \bar{b}_0\bar{b}_1$.
        \item[(iii)] For every clique $i_0<...<i_{k-1}$ in $G$, $b_{i_0}^0b_{i_0}^1...b_{i_{k-1}}^0b_{i_{k-1}}^1$ is $\ind$-independent over $C$. 
        \end{enumerate}
        For the base case, take $b_0^0 = \bar{b}_2$. Since $\bar{b}_1$ is algebraically closed, by full existence, we can find $b_0^1 \equiv_{\bar{b}_1} \bar{b}_0$ such that $b_0^1 \ind_{\bar{b}_1}b_0^0$. Since $\bar{b}_1 \cap b_0^1b_0^0 = (\bar{b}_1 \cap b_0^1)(\bar{b}_1\cap b_0^0) = C(\bar{b}_1 \cap \bar{b}_2) = C$, by freedom we have $b_0^1 \ind_C b_0^0.$ Thus, condition (iii) holds. Moreover, condition (i) says $b_0^0 \ind_C \emptyset$ which follows from $b_0^1 \ind_Cb_0^0$ by symmetry and monotonicity, while condition (ii) is vacuously true. Notice that, by Lemma \ref{free-am-implies-alg-ind}, since $C \subseteq b_0^0\cap b_0^1$ and $b_0^1 \ind_C b_0^0$, we have that $C = b_0^0 \cap b_0^1$. 
        Assume now that for a graph $G = (n,E)$ we have $(b_i^0,b_i^1)_{i<n}$ satisfying (i),(ii) and (iii). Let $A,B \subseteq n$, $A\cap B = \emptyset$, $A\cup B = n$. Define $G(A,B) = (n+1,E)$ an extension of $G$ where $n E i$ for all $i \in B$ while $\neg n E j$ for all $j \in A$. By full existence, find $b_n^0 \equiv_C b_0^0 = \bar{b}_2$ such that $b_n^0 \ind_C b_{<n}^0b_{<n}^1$. Thus condition (i) is met. Since $(b_i^0)_{i \in A}$ is $\ind$-independent over $C$, by Lemma \ref{dividingalongIndepSequences}, there is $b^*$ such that 
        \[b_i^0b^* \equiv_C \bar{b}_0\bar{b}_1 \quad \forall i \in A\]
        By the axiom of closure, $(b_i^0)_{i \in A}$ is closed and thus,
        by full existence, find $b_n^1 \equiv_{(b_i^0)_{i \in A}}b^*$ such that \[b_n^1 \ind_{(b_i^0)_{i \in A}} b_n^0b_{<n}^1(b_i^0)_{i \in B}\] 
        Notice that, for $i<j<n$ and $\epsilon,\delta \in 2$, for $b_i^\epsilon$ and $b_j^\delta$ we either have $b_i^\epsilon b_j^\delta \equiv_C \bar{b}_0\bar{b}_1$ or $b_i^\epsilon \ind_C b_j^\delta$; in either case, by Lemma \ref{free-am-implies-alg-ind}, we have $b_i^\epsilon \cap b_j^\delta = C$. Similarly, for $i \in A$, $b_n^1 \cap b_i^0 = C$ and $b_n^0 \cap b_i^0 = C$.   
        Thus, by freedom, $b_n^1 \ind_C b_n^0b_{<n}^1(b_i^0)_{i \in B}$. This, the condition on $b^*$ and inductive hypothesis imply both (ii) and (iii).
        
        Finally, we prove that $\fhi(x;y)$ is $\PM^{(2)}$. Suppose $G = (n,E)$ is any finite graph. By what we've done so far, we can find $(b_i^0b_i^1)_{i<n}$ satisfying (i),(ii) and (iii). Since each $b_i^\epsilon$ is a copy of the algebraic closure of an element of the sequence $(b_i)_{i<\omega}$ witnessing $\SOP_3$, call $c_i^\epsilon$ the subtuple of $b_i^\epsilon$ corresponding to $b_i$.
        We claim that $(c_i^0c_i^1)_{i<\omega}$ realizes $G$ with $\fhi(x;y)$. Suppose $i<j<n$ are such that $\neg i E j $ in $G$. Then, by (ii), $b_i^0b_j^1 \equiv_C\bar{b}_0\bar{b}_1$ and thus $c_i^0c_j^1 \equiv b_0b_1$. In particular, $\psi(c_i^0,c_j^1)$ holds. So, if $\fhi(x;c_i^0c_i^1) \land \fhi(x;c_j^0c_j^1)$ was realized by some $a$, we would get $\psi(a,c_i^0)\land \psi(c_i^0,c_j^1) \land \psi(c_j^1,a)$ contradicting $\SOP_3$. 
        
        Now, assume $i_0<...<i_{k-1}$ is a clique in $G$. By (iii), $b_{i_0}^0b_{i_0}^1...b_{i_{k-1}}^0b_{i_{k-1}}^1$ is $\ind$-independent over $C$. By Lemma \ref{pairsofIndeptuples}, the sequence $(b_{i_l}^0b_{i_l}^1)_{l<k}$ is $\ind$-independent over $C$. 
        Since $b_{i_0}^0 \equiv_C b_0^0$, there is $\sigma \in \Aut(\U/C)$ such that $\sigma(b_{i_0}^0) = b_0^0$. Thus $b_{i_0}^0b_{i_0}^1 \equiv_C b_0^0\sigma(b_{i_0}^1)$. Now, since $b_{i_0}^1 \ind_C b_{i_0}^0$, by invariance $\sigma(b_{i_0}^1) \ind_C b_0^0$; moreover, $b_0^1 \ind_Cb_0^0$; thus, by stationarity $b_{i_0}^0b_{i_0}^1 \equiv_C b_0^0 \sigma(b_{i_0}^1) \equiv_C b_0^0b_0^1 = \bar{b}_2b_0^1$. Moreover $b_0^1 \equiv_{\bar{b}_1} \bar{b}_0$. Thus $\psi(b_1,b_0^0)\land \psi(b_0^1,b_1)$. So we can find an algebraically closed tuple $a$ with $a b_{i_0}^0b_{i_0}^1 \equiv_C \bar{b}_1 b_0^0b_0^1$ and so $\fhi(a;b_{i_0}^0b_{i_0}^1)$. Moreover, $a \cap b_{i_0}^0b_{i_0}^1 = C$. Thus, by Lemma \ref{dividingalongIndepSequences}, there is $a^*$ such that $a^*b_{i_l}^0b_{i_l}^1 \equiv_C ab_{i_0}^0b_{i_1}^0$ for all $l<k$. In particular 
        \[\models\bigwedge_{l<k}\fhi(a^*,b_{i_l}^0b_{i_l}^1)\]
    
    In conclusion $\fhi(x;y)$ is $\PM^{(2)}$. 
\end{proof}

\begin{defn}
    A theory is \emph{modular} if $\ind^a$ satisfies base monotonicity, where 
    \begin{itemize}
        \item $A \ind^a_C B \iff \acl(AC)\cap \acl(BC) = \acl(C)$
        \item a ternary relation $\ind$ satisfies \emph{base monotonicity} if, for all $A,B,C,D$ such that $D \subseteq C \subseteq B$, 
        \[A \ind_D B \implies A \ind_C B \]
    \end{itemize}
    Alternatively, $T$ is modular if the lattice of algebraically closed sets is modular. 
\end{defn}

In \cite{free-amalgamation_Conant}, Conant proved that a modular free amalgamation theory is simple if and only if it is $\NSOP_3$. Using this, together with our previous theorem, we get
\begin{cor}
    Let $T$ be a modular free amalgamation theory. Then $T$ is simple if and only if $T$ is $\NPM^{(2)}$.
\end{cor}
\begin{proof}
    Clearly, if $T$ is simple, then $T$ is $\NPM^{(2)}$. Indeed, since $2$-$\TP$ is described by positive patterns with inconsistencies of size $2$, it is implied by $\PM^{(2)}$.

    On the other hand, if $T$ is non-simple, by Theorem 1.17 from \cite{free-amalgamation_Conant}, $T$ has $\SOP_3$. Then, by Theorem \ref{freeamwithSOP3isPM^2}, $T$ is $\PM^{(2)}$.  
\end{proof}

Given the fact that $\NSOP_1   =\NSOP_3$ (see \cite{OnNSOP2theories},\cite{chernikov2026sop2equalssop3}), using Theorem \ref{freeamwithSOP3isPM^2}, we trivially have:
\begin{cor}\label{trivial-corollary-of-dichotomy-for-freeAm}
    Let $T$ be a free amalgamation theory. Then $T$ is $\NSOP_1$ if and only if $T$ is  $\NPM^{(2)}$. 
\end{cor}

Before the Chernikov's proof of $\NSOP_2 = \NSOP_3$ (\cite{chernikov2026sop2equalssop3}),  Mutchnik showed that free amalgamation theories are either $\NSOP_1$ or $\SOP_3$ (see \cite[Theorem 4.4]{mutchnik2024conantindependencegeneralizedfreeamalgamation}). This result, together with Theorem \ref{freeamwithSOP3isPM^2}, gave a proof of Corollary \ref{trivial-corollary-of-dichotomy-for-freeAm} even before we knew that $\NSOP_3 = \NSOP_1$.

\begin{proof}
    If $T$ is $\NSOP_1$, then it is clearly $\NPM^{(2)}$. On the other hand, if $T$ has $\SOP_1$, then by Theorem 4.4 in \cite{mutchnik2024conantindependencegeneralizedfreeamalgamation}, it has $\SOP_3$ and so, by Theorem \ref{freeamwithSOP3isPM^2}, it has $\PM^{(2)}$.
\end{proof}

%Other-partial-results
\section{Other partial results}
\label{other-partial-results}

In this section we prove other results providing some evidence for a positive answer to the dichotomy question (Question \ref{question}). 

\subsection{Using the real-valued strict-order property hierarchy 
}

In \cite{mutchnik2025approximationsstrictorderproperty}, Mutchnik defined a possible extension of the $\NSOP_n$ hierarchy where $n$ can be chosen to be any real number $r \geq 3$, as follows:
\begin{defn}[\cite{mutchnik2025approximationsstrictorderproperty}]
Fix $r \geq 3$. A theory $T$ is $\NSOP_r$ if there is no definable relation $R(x,y)$ such that:
\begin{itemize}
    \item there exists a sequence $(a_i)_{i<\omega}$ such that $R(a_i,a_j)$ holds for every $i<j<\omega$, but
    \item there is no set $(a_\theta)_{\theta \in S^1}$, indexed by the unit circle $S^1$, such that $R(a_\theta,a_\psi)$ hold for all $\theta,\psi \in S^1$ with $\psi$ lying at most $\dfrac{2\pi}{r}$ radians counterclockwise from $\theta$.
\end{itemize}
Otherwise, we say $T$ has $\SOP_r$. 
\end{defn}
While in \cite{mutchnik2025approximationsstrictorderproperty}, Mutchnik focused on the question of whether the real-valued $\NSOP_r$ hierarchy properly extend the integer-valued $\NSOP_n$ hierarchy, in \cite{mutchnik2026applicationsrealstrictorder} they present a general way in which this could fail: 

\begin{hyp}[See \cite{mutchnik2026applicationsrealstrictorder}]
    Fix $N \geq 3$. \emph{The real-valued $\NSOP_r$ hierarchy is non-distinct from the integer-valued $\NSOP_n$ hierarchy between $\NSOP_N$ and $\NSOP_{N+1}$ on sufficiently general grounds} (shortly denoted by $\mathbf{H}(N)$) if for all finite directed graphs $G'$ without directed cycles of length $\leq N$ and for all infinite, directed graphs $G$ containing some sequence $(v_i)_{i<\omega}$ with $v_i E v_j$ for all $i<j<\omega$ and with no quantifier-free formula in the graph language exhibiting $\SOP_{N+1}$, there is a weak embedding (i.e.\ an injective homomorphism) $G' \to G$ . 

    We say that \emph{the real-valued $\NSOP_r$ hierarchy is non-distinct from the integer-valued $\NSOP_n$ hierarchy on sufficiently general grounds} (shortly denoted $\mathbf{H}$) if $H(n)$ holds for every $n\geq 3$. 
\end{hyp}
Again, in \cite{mutchnik2026applicationsrealstrictorder} (see Theorem 4.4), they proceed to prove that, assuming $\mathbf{H}$, we get the collapse of the $\NSOP_n$ hierarchy inside $\NTP_2$. That is, for every $n\geq 3$, 
\[\NSOP_n \cap \NTP_2 = \NSOP_{n+1}\cap \NTP_2\]

Following a similar approach, we can, in a sense, approximate a positive answer to the dichotomy question as close as we want and collapse the $\NSOP_n$ hierarchy inside $\NPM^{(2)}$. More precisely:  

\begin{thm}
    For every $2<n<\omega$, assuming $\mathbf{H}(n)$, we have that 
    \[\NSOP_n \cap \NPM^{(2)} = \NSOP_{n+1}\cap \NPM^{(2)}\]
    Furthermore, assuming $\mathbf{H}$, we have that for all $2<n<\omega$
    \[\NSOP_3 = \NSOP_n \cap \NPM^{(2)}\]
\end{thm}

\begin{proof}
    Fix $N \geq 3$ and assume $\mathbf{H}(N)$. Since $\NSOP_N \subseteq \NSOP_{N+1}$, we know that $\NSOP_N \cap \NPM^{(2)} \subseteq \NSOP_{N+1}\cap \NPM^{(2)}$. We want to show the reverse inclusion. Suppose, $T \in \NSOP_{N+1}\cap \NPM^{(2)}$ but, toward contradiction, it has $\SOP_N$. We will get a contradiction by showing that $T$ has a $\PM^{(2)}$ formula. 

    Let $R(x,y)$ the (without loss of generality antisymmetric and irreflexive) formula witnessing $\SOP_N$ with $(a_i)_{i<\omega}$. So, for all $i<j<\omega$, we have $R(a_i,a_j)$ and 
    \[\neg \exists x_0...x_{N-1}\bigwedge_{i<N-1}R(x_i,x_{i+1}) \land R(x_{N-1},x_0)\]
    We claim that the formula
    \[\fhi_R(\bar{x},yz) = \bigwedge_{i<N-2}R(x_i,x_{i+1})\land R(y,x_0)\land R(x_{N-1},z)\]
    is $\PM^{(2)}$. 

    Fix a finite graph $G = (V,E)$. Define a directed graph $G' = (V',E')$ as follows
    \[V' = \{v_i' \: i \in 2, v \in V\}\cup\{a_C^i \: i < N-2,\ C \text{ clique in $G$}\}\]
    with directed edges
    \[v_i'E'w_j' \iff v\neq w \land i=1 \land j=0\land \neg vEw\]
    and, for all cliques $C$ in $G$
    \[\bigwedge_{v \in C}\fhi_{E'}(\bar{a}_C,v_0'v_1')\]
    where $\bar{a}_C = a_C^0...a_C^{N-1}$, and no other directed edges. 

    We now show that $G'$ has no directed cycle of length $\leq N$. 

    Suppose $O$ is a directed cycle in $G'$. Since $\{v_i' \: v \in V, i\in 2\}$ is bipartite and the directed edges in it are always from vertices of index $1$ to vertices of index $0$, we know that $O \not \subseteq \{v_i' \: v \in V, i \in 2\}$. Thus, it must involve some of the $\bar{a}_C$'s. But, for a fixed $C$, the tuples in $\bar{a}_C$ are only connected to each others (forming a directed path) and to the independent set $\{v_i' \: v \in C\}$. Thus multiple $\bar{a}_C$ must be involved in $O$ and so the length of $O$ must exceed $N$. In conclusion $G'$ has no directed cycle of length $\leq N$. 

    Now, by $\mathbf{H}(N)$ and by the assumption that $T$ is $\NSOP_{N+1}$, we can (weakly) embed (maintaining the same names for the vertices) $G'$ into the directed graph $(\U,R)$ (where $\U \models T$). Given $v \neq w$ in $G$ such that $\neg v E w$, we have that $v_1'Rw_0'$. Thus, if there was $\bar{a}$ such that 
    \[\models \fhi_R(\bar{a},v_0'v_1') \land \fhi_R(\bar{a},w_0'w_1')\]
    we would have that 
    \[v_1'\bar{a}w_0' \models \bigwedge_{i<N-1}R(x_i,x_{i+1})\land R(x_{N-1},x_0)\]
    which is inconsistent by $\SOP_N$. 
    On the other hand, if $C$ is a clique in $G$, then 
    \[\bar{a}_C \models \bigwedge_{v \in C}\fhi_R(\bar{x},v_0'v_1')\]
    In conclusion $\fhi_R$ would be $\PM^{(2)}$, contradicting the fact that $T$ is $\NPM^{(2)}$. Thus $T$ is $\NSOP_{N}$.

    Furthermore, since $\NSOP_3 = \NSOP_3 \cap \NPM^{(2)}$, assuming $\mathbf{H}$ and inductively use what we have just proved, we have that for all $2<n<\omega$
    \[\NSOP_3 =  \NSOP_n \cap \NPM^{(2)}\]
\end{proof}

\subsection{$\SOP$ formula from omitting specific graphs}
First, there is a characterization of interval graph by forbidding a set of graphs: 

\begin{prop}[\cite{Interval-graph_forbidden-subgraphs}, Theorem 4]\label{forbidden-subgraph-for-interval-thm}
A finite graph is interval if and only if it does not contain any cycle with $n\geq 4$ vertices and any $T_2$, $X_{31}$, $(XF_2^m)_{m\geq 1}$ and $(XF_3^n)_{n \geq1}$ (shown, in this order, below).
        \begin{center}
         \begin{tikzpicture}[node distance={8mm}, main/.style = {draw, circle}] 
    \node[main] (1){};
    \node[main] (2) [right of=1]{};
    \node[main] (3) [right of=2]{};
    \node[main] (4) [right of=3]{};
    \node[main] (5) [right of=4]{};
    \node[main] (6) [below of=3]{};
    \node[main] (7) [below of=6]{};
    \draw (1)--(2);
    \draw (2)--(3);
    \draw (3)--(4);
    \draw (4)--(5);
    \draw (6)--(3);
    \draw (7)--(6);
    \end{tikzpicture} 
    \end{center}
    \begin{center}
    \rotatebox[]{180}{ \begin{tikzpicture}[node distance={8mm}, main/.style = {draw, circle}] 
    \node[main] (1){};
    \node[main] (2) [below of=1]{};
    \node[main] (3) [right of=2]{};
    \node[main] (4) [left of=2]{};
    \node[main] (5) [below of=2]{};
    \node[main] (6) [left of=5]{};
    \node[main] (7) [right of=5]{};
    \draw (1)--(2);
    \draw (3)--(2);
    \draw (4)--(2);
    \draw (3)--(7);
    \draw (4)--(6);
    \draw (5)--(6);
    \draw (5)--(7);
    \draw (3)--(5);
    \draw (4)--(5);
    \draw (2)--(5);
    \end{tikzpicture}}
    \end{center}

    \begin{center}
     \begin{tikzpicture}[node distance={10mm}, main/.style = {draw, circle}] 
    \node[main] (1){};
    \node[main] (2) [right of=1]{};
    \node[] (3) [right of=2]{};
    \node[] (middle) [right of=3]{};
    \node[] (3a) [right of =middle]{};
    \node[main] (4) [right of=3a]{};
    \node[main] (5) [right  of=4]{};
    \node[main] (6) [below of=middle]{};
    \node[main] (7) [below of=6]{};
    \draw[dashed] (2)--(4) node[midway, font=\scriptsize, sloped, above] {$m+1$ many};
    \draw (1)--(2);
    \draw (4)--(5);
    \draw (2)--(6)--(7);
    \draw (4)--(6);
    \draw (6)-- (3);
    \draw (6)--(3a);
    \draw (6)--(middle);
    \end{tikzpicture}
    \end{center}

    \begin{center}
        \begin{tikzpicture}[node distance={10mm}, main/.style = {draw, circle}]
            \node[main] (1){};
            \node[main] (2)[below of =1,right of =1]{};
            \node[main] (3) [below of =1,left of =1]{};
            \node[main] (4) [below of =3, left of =3]{};
            \node[main] (5)[below of =2,right of =2]{};
            \node[main] (6) [below of =3]{};

            \node[main] (8) [left of =5]{};
            \draw (1)--(2)--(5)--(8);
            \draw (1)--(3)--(4)--(6)--(3)--(2)--(8);
            \draw (3)--(8);
            \draw (2)--(6);
        
             \draw[dashed] (6)--(8) node[midway, font=\scriptsize, sloped, below] {$n+1$ many};

            \node[] (phantom) [left of=8]{};
            \draw (2)--(phantom)--(3);
        \end{tikzpicture}
    \end{center}
\end{prop}

\begin{rem}
Assume that we have a theory $T$ which is $\SOP_3$, witnessed by a formula $\fhi(x;y)$ that realizes every finite interval graph (see Remark \ref{rem-SOP3=real-Intervals}) but $T$ is $\NPM^{(2)}$. In particular, for some finite (non-interval) graph $G$, $\fhi$ doesn't realize $G$. If in this situation we can prove the existence of a formula with $\SOP$, we get a positive answer to Question \ref{question}. 
\end{rem}

The following results show how to get an $\SOP$ formula assuming that a formula $\fhi(x;y)$ realizes every interval graph but omits one of the graphs described in the previous proposition. 

In the proofs of these results, starting from a graph that the formula doesn't realize, we are going to build an interval graph with the idea to capture a strict chain of definable sets (which witnesses $\SOP$). 

For instance, in the proof of the next proposition, we build a graph made of rows $(a_ib_ic_id_ie_i)$ which form a path; then, for $i<j$, we connect every element of the $i$th row with the central element of the $j$th row (i.e.\ $c_j$). 

For instance, looking at the $i$th and $j$th rows, for $i<j$, we would have
\begin{center}
    \begin{tikzpicture}[node distance={15mm}, main/.style = {draw, circle}] 
    \node[] (1){$a_j$};
    \node[] (2) [right of=1]{$b_j$};
    \node[] (3) [right of=2]{$c_j$};
    \node[] (4) [right of=3]{$d_j$};
    \node[] (5) [right of=4]{$e_j$};
    \node[] (6) [below of=1]{$a_i$};
    \node[] (7) [below of=2]{$b_i$};
    \node[] (8) [below of=3]{$c_i$};
    \node[] (9) [below of=4]{$d_i$};
    \node[] (10) [below of=5]{$e_i$};
   
    \draw (1)--(2);
    \draw (2)--(3);
    \draw (3)--(4);
    \draw (4)--(5);
    \draw (6)--(3);
    \draw (7)--(6);
    \draw (3)--(7);
    \draw (3)--(8);
    \draw (3)--(9);
    \draw (3)--(10);
    \draw (8)--(7);
    \draw (8)--(9);
    \draw (9)--(10);

    \end{tikzpicture} 
    \end{center}

Interpreting each vertex as an interval on the line (say $x_i = [x_i^l,x_i^r]$ for $x \in \{a,b,c,d,e\}$), each row will look like 
\[a_i^l < b_i^l < a_i^r < c_i^l < b_i^r < d_i^l < c_i^r< e_i^l<d_i^r<e_i^r \]
Then, the previous row will be all contained in the open interval $(b_i^r,d_i^l)$ which is fully contained in $c_i = [c_i^l,c_i^r]$.
This shows, first of all, that the graph is interval; secondly, any non-empty intersection with some of the intervals in the $(i-1)$th row, will imply a non-empty intersection with some of the intervals in the $i$th row, hence the inclusion of definable subsets.

\begin{prop}
        Suppose $\fhi(x;y)$ realizes every finite interval graph but omtis $T_2$ (shown below)
        \begin{center}
         \begin{tikzpicture}[node distance={8mm}, main/.style = {draw, circle}] 
    \node[main] (1){};
    \node[main] (2) [right of=1]{};
    \node[main] (3) [right of=2]{};
    \node[main] (4) [right of=3]{};
    \node[main] (5) [right of=4]{};
    \node[main] (6) [below of=3]{};
    \node[main] (7) [below of=6]{};
    \draw (1)--(2);
    \draw (2)--(3);
    \draw (3)--(4);
    \draw (4)--(5);
    \draw (6)--(3);
    \draw (7)--(6);
    \end{tikzpicture} 
    \end{center}
    Then, there is a formula with $\SOP$. 
    \end{prop}
    \begin{proof}
         Define the graph $G = (V,E)$ with vertices $(a_i,b_i,c_i,d_i,e_i)_{i<n}$ where, for each $i<n$, $\{a_i,b_i\},\{b_i,c_i\},\{c_i,d_i\},\{d_i,e_i\} \in E$ and, for each $j<i$, $\{c_i,x_j\} \in E$ for $x \in \{a,b,c,d,e\}$.
       For example, for $n = 2$ we have
       \begin{center}
    \begin{tikzpicture}[node distance={15mm}, main/.style = {draw, circle}] 
    \node[] (1){$a_1$};
    \node[] (2) [right of=1]{$b_1$};
    \node[] (3) [right of=2]{$c_1$};
    \node[] (4) [right of=3]{$d_1$};
    \node[] (5) [right of=4]{$e_1$};
    \node[] (6) [below of=1]{$a_0$};
    \node[] (7) [below of=2]{$b_0$};
    \node[] (8) [below of=3]{$c_0$};
    \node[] (9) [below of=4]{$d_0$};
    \node[] (10) [below of=5]{$e_0$};
   
    \draw (1)--(2);
    \draw (2)--(3);
    \draw (3)--(4);
    \draw (4)--(5);
    \draw (6)--(3);
    \draw (7)--(6);
    \draw (3)--(7);
    \draw (3)--(8);
    \draw (3)--(9);
    \draw (3)--(10);
    \draw (8)--(7);
    \draw (8)--(9);
    \draw (9)--(10);

    \end{tikzpicture} 
    \end{center}
    For any $n<\omega$, this is an interval graph so $\fhi(x;y)$ realizes it (let's maintain the same name of the vertices for the realizations). 
    Consider the edge relation
    \[\psi(y_0,y_1) = \exists x \fhi(x;y_0) \land \fhi(x;y_1) \land y_0 \neq y_1\]
    and define the formula 
    \[\theta(x;y_0y_1y_2y_3y_4) = \bigvee_{i<5}\psi(x,y_i)\]
    Fix $i<j<n$. Notice that \[a_j \models \psi(x,b_j) \land \bigwedge_{f \in \{a,b,c,d,e\}}\neg \psi(x,f_i)\]
    Thus \[a_j\models \theta(x;a_jb_jc_jd_je_j)\land \neg \theta(x;a_ib_ic_id_ie_i)\]
    Now, if $a \models \theta(x;a_ib_ic_id_ie_i)$ then, for some $f \in \{a,b,c,d,e\}$, $a \models \psi(a,f_i)$. But now $a_jb_jc_jd_je_jf_ia$ form a copy of $T_2$ unless there is an edge from $a$ to $f_j$ for $f \in \{a,b,c,d,e\}$. 
    Thus $a \models \theta(x;a_jb_jc_jd_je_j)$. In conclusion, for every $i<j<n$
    \[\theta(\U;a_ib_ic_id_ie_i) \subsetneq \theta(\U;a_jb_jc_jd_je_j)\]
    By compactness, $\theta$ has $\SOP$. 
    \end{proof}

    \begin{prop}
         Suppose $\fhi(x,y)$ realizes every finite interval graph but omits $X_{31}$ (shown below)
    \begin{center}
    \rotatebox[]{180}{ \begin{tikzpicture}[node distance={8mm}, main/.style = {draw, circle}] 
    \node[main] (1){};
    \node[main] (2) [below of=1]{};
    \node[main] (3) [right of=2]{};
    \node[main] (4) [left of=2]{};
    \node[main] (5) [below of=2]{};
    \node[main] (6) [left of=5]{};
    \node[main] (7) [right of=5]{};
    \draw (1)--(2);
    \draw (3)--(2);
    \draw (4)--(2);
    \draw (3)--(7);
    \draw (4)--(6);
    \draw (5)--(6);
    \draw (5)--(7);
    \draw (3)--(5);
    \draw (4)--(5);
    \draw (2)--(5);
    \end{tikzpicture}}
    \end{center}
    Then there exists a formula with $\SOP$. 
    \end{prop}

    \begin{proof}
        Fix $n<\omega$ and define a graph $G$ with vertices $(a_ib_ic_id_ie_i)_{i<n}$ such that, for each $i<n$, $a_iEb_iEc_iEd_iEe_i$ and $a_iEc_iEe_i$. Moreover, for $j<i$, we have $x_jEb_i \land x_j Ec_i \land x_jEd_i$ for $x \in \{a,b,c,d,e\}$.
        For example, for $n=2$, we have
        \begin{center}
        \begin{tikzpicture}[node distance={15mm}]
            \node[] (a1){$a_1$};
            \node[] (b1) [right of =a1] {$b_1$};
            \node[] (c1) [right of =b1] {$c_1$};
            \node[] (d1) [right of =c1] {$d_1$};
            \node[] (e1) [right of =d1] {$e_1$};

            \node[] (a0)[below of=a1]{$a_0$};
             \node[] (b0) [right of =a0] {$b_0$};
            \node[] (c0) [right of =b0] {$c_0$};
            \node[] (d0) [right of =c0] {$d_0$};
            \node[] (e0) [right of =d0] {$e_0$};

            \draw (a1)--(b1)--(c1)--(d1)--(e1);
            \draw (a0)--(b0)--(c0)--(d0)--(e0);

            \draw(a1) to [out=50](c1);
            \draw(c1) to [out=50](e1);
            \draw(a0) to [out=-50,in=-120](c0);
            \draw(c0) to [out=-50,in=-120](e0);

            \draw (a0)--(b1);
            \draw (a0)--(c1);
            \draw (a0)--(d1);

            \draw (b0)--(b1);
            \draw (b0)--(c1);
            \draw (b0)--(d1);

            \draw (c0)--(b1);
            \draw (c0)--(c1);
            \draw (c0)--(d1);
            
            \draw (d0)--(b1);
            \draw (d0)--(c1);
            \draw (d0)--(d1);

            \draw (e0)--(b1);
            \draw (e0)--(c1);
            \draw (e0)--(d1);
        \end{tikzpicture}
        \end{center}

        For any $n<\omega$, this is an interval graph and so $\fhi(x;y)$ realizes it (let's, again, maintain the same name of the vertices for the realizations). Consider the edge relation 
        \[\psi(y_0,y_1) = \exists x \fhi(x;y_0)\land \fhi(x;y_1) \land y_0 \neq y_1\]
        and define the formula

    \[\theta(x;\bar{y}) = \bigvee_{i<4}\psi(x,y_i)\]
    Fix $i<j<n$. Notice that \[a_j \models \psi(x,b_j) \land \bigwedge_{f \in \{a,b,c,d,e\}}\neg \psi(x,f_i)\] Thus \[a_j \models \theta(x;a_jb_jc_jd_je_j)\land \neg \theta(x;a_ib_ic_id_ie_i)\]
    Now, if $a \models \theta(x;a_ib_ic_id_ie_i)$ then, for some $f \in \{a,b,c,d,e\}$, $a \models \psi(a,f_i)$. But now $\fhi$ realizes $X_{31}$, as witnessed by $a,f_i,a_j,b_j,c_j,d_j,e_j$. For instance, if $f = c$, we get

    \begin{center}
        \begin{tikzpicture}[node distance={15mm}]
            \node[] (aj){$a_j$};
            \node[] (bj) [right of =aj] {$b_j$};
            \node[] (cj) [right of =bj] {$c_j$};
            \node[] (dj) [right of =cj] {$d_j$};
            \node[] (ej) [right of =dj] {$e_j$};

            \node[] (ai)[below of=aj]{$a_i$};
             \node[] (bi) [right of =ai] {$b_i$};
            \node[] (ci) [right of =bi] {$c_i$};
            \node[] (di) [right of =ci] {$d_i$};
            \node[] (ei) [right of =di] {$e_i$};

            \draw[thick,red] (aj)--(bj)--(cj)--(dj)--(ej);
            \draw (ai)--(bi)--(ci)--(di)--(ei);

            \draw[thick,red](aj) to [out=50](cj);
            \draw[thick,red] (cj) to [out=50](ej);
            \draw(ai) to [out=-50,in=-120](ci);
            \draw(ci) to [out=-50,in=-120](ei);

            \draw[thick,red] (ci)--(bj);
            \draw[thick,red] (ci)--(cj);
            \draw[thick,red] (ci)--(dj);

            \node[] (a) [below of =ci]{$a$};
            \draw[thick,red] (a)--(ci);
        \end{tikzpicture}
        \end{center}
        So there must be an edge between $a$ and $x_j$ for some $x \in \{a,b,c,d,e\}$. Thus $a \models \theta(x;a_jb_jc_jd_je_j)$. 
        In conclusion, for every $i<j<n$, we have 
        \[\theta(\U;a_ib_ic_id_ie_i) \subsetneq \theta(\U;a_jb_jc_jd_je_j)\]
    By compactness, $\theta$ has $\SOP$.    
    \end{proof}

    \begin{prop}
        Suppose $\fhi(x;y)$ realizes every finite interval graph but, for some $m\geq 1$ omits $XF^m_2$ (shown below)
        \begin{center}
     \begin{tikzpicture}[node distance={10mm}, main/.style = {draw, circle}] 
    \node[main] (1){};
    \node[main] (2) [right of=1]{};
    \node[] (3) [right of=2]{};
    \node[] (middle) [right of=3]{};
    \node[] (3a) [right of =middle]{};
    \node[main] (4) [right of=3a]{};
    \node[main] (5) [right  of=4]{};
    \node[main] (6) [below of=middle]{};
    \node[main] (7) [below of=6]{};
    \draw[dashed] (2)--(4) node[midway, font=\scriptsize, sloped, above] {$m+1$ many};
    \draw (1)--(2);
    \draw (4)--(5);
    \draw (2)--(6)--(7);
    \draw (4)--(6);
    \draw (6)-- (3);
    \draw (6)--(3a);
    \draw (6)--(middle);
    \end{tikzpicture}
    \end{center}

    Then, there is a formula with $\SOP$. 
    \end{prop}
    \begin{proof}
        For simplicity, we are going to prove it in the case $m=1$. The proof for $m>1$ is similar.
        So, we are assuming $\fhi(x;y)$ realizes every interval graph but doesn't realize $XF^1_2$
        \begin{center}
     \begin{tikzpicture}[node distance={8mm}, main/.style = {draw, circle}] 
    \node[main] (1){};
    \node[main] (2) [right of=1]{};
    \node[] (3) [right of=2]{};
    \node[main] (4) [right of=3]{};
    \node[main] (5) [right  of=4]{};
    \node[main] (6) [below of=3]{};
    \node[main] (7) [below of=6]{};
    \draw (1)--(2)--(4)--(5);
    \draw (2)--(6)--(7);
    \draw (4)--(6);
    \end{tikzpicture}
    \end{center}

    Fix $n<\omega$ and define a graph $G$ with vertices $(a_ib_ic_id_i)_{i<n}$ such that, for each $i<n$, $a_i E b_i E c_i E d_i$ and, for $i<j$, $x_i E b_j$ and $x_i E c_j$ for all $x \in \{a,b,c,d\}$. For example, for $n=2$, we have
    \begin{center}
        \begin{tikzpicture}[node distance={15mm}]
        \node[] (aj){$a_1$};
        \node[] (bj)[right of =aj]{$b_1$};
        \node[] (cj)[right of =bj]{$c_1$};
        \node[] (dj)[right of =cj]{$d_1$};

        \node[] (ai)[below of =aj]{$a_0$};
        \node[] (bi)[right of =ai]{$b_0$};
        \node[] (ci)[right of =bi]{$c_0$};
        \node[] (di)[right of =ci]{$d_0$};

        \draw (aj)--(bj)--(cj)--(dj);
        \draw (ai)--(bi)--(ci)--(di);

        \draw (ai)--(bj);
        \draw (ai)--(cj);

        \draw (bi)--(bj);
        \draw (bi)--(cj);

        \draw (ci)--(bj);
        \draw (ci)--(cj);

        \draw (di)--(bj);
        \draw (di)--(cj);

        \end{tikzpicture}
    \end{center}
    For any $n<\omega$, this is an interval graph and so $\fhi(x;y)$ realizes it (let's, once again, maintain the same name of the vertices for the realizations). Consider the edge relation 
        \[\psi(y_0,y_1) = \exists x \fhi(x;y_0)\land \fhi(x;y_1) \land y_0 \neq y_1\]
        and define the formula

    \[\theta(x;\bar{y}) = \bigvee_{i<3}\psi(x,y_i)\]

    Fix $i<j<n$. Notice that \[a_j \models \psi(x,b_j) \land \bigwedge_{f \in \{a,b,c,d\}}\neg \psi(x,f_i)\] Thus \[a_j \models \theta(x;a_jb_jc_jd_j)\land \neg \theta(x;a_ib_ic_id_i)\]
    Now, if $a \models \theta(x;a_ib_ic_id_i)$ then, for some $f \in \{a,b,c,d\}$, $a \models \psi(a,f_i)$. But now $\fhi$ realizes $XF_{2}^1$, as witnessed by $a,f_i,a_j,b_j,c_j,d_j$. For instance, if $f = b$, we get
    \begin{center}
        \begin{tikzpicture}[node distance={15mm}]
        \node[] (aj){$a_j$};
        \node[] (bj)[right of =aj]{$b_j$};
        \node[] (cj)[right of =bj]{$c_j$};
        \node[] (dj)[right of =cj]{$d_j$};

        \node[] (ai)[below of =aj]{$a_i$};
        \node[] (bi)[right of =ai]{$b_i$};
        \node[] (ci)[right of =bi]{$c_i$};
        \node[] (di)[right of =ci]{$d_i$};

        \draw[thick,red] (aj)--(bj)--(cj)--(dj);
        \draw (ai)--(bi)--(ci)--(di);

        \draw[thick,red] (bi)--(bj);
        \draw[thick,red] (bi)--(cj);

        \node[] (a) [below of =bi]{$a$};
        \draw[thick,red] (a)--(bi);

        \end{tikzpicture}
    \end{center}
    So there must be an edge between $a$ and $x_j$ for some $x \in \{a,b,c,d\}$. Thus $a \models \theta(x;a_jb_jc_jd_j)$. In conclusion
    \[\theta(\U;a_ib_ic_id_i)\subsetneq \theta(\U;a_jb_jc_jd_j)\]
    By compactness, $\theta$ has $\SOP$.
    
    \end{proof}
\begin{rem}
A similar approach doesn't seem to work for $n$-cycles ($n\geq 4$) and for $(XF_3^n)_{n\geq 1}$. In particular, for each $n\geq 1$, $XF_3^n$ has no vertex of degree $1$ and so a formula $\theta(x,\bar{y})$ that impose a realization of $XF_3^n$ on a row of an interval graph, would require a consistency of at least $3$ vertices while implying the consistency of only two.  

We also remark that, even if a similar method worked for cycles and $XF_3^n$, we still wouldn't have solved the problem. Indeed, not realizing a graph containing one of the forbidden subgraphs doesn't \emph{a priori} imply not realizing such subgraph. 
\end{rem}

%Acknowledgments
\addtocontents{toc}{\protect\setcounter{tocdepth}{0}}
\section*{Acknowledgements}
Most of this work is part of the author's PhD dissertation (\cite{MicheleBailettiPhDthesis}). The author sincerely thanks his former advisor Alex Kruckman, for productive discussions, enlightening comments and suggestions. Furthermore, this work was partially supported by Project 24-12591M of the Czech Science Foundation (GAČR), and by the long-term strategic development financing of the Institute of Computer Science (RVO: 67985807).

%AI Statement
\section*{AI Statement} No AI tool was used in this paper, neither in the discovery nor in the writing process.

%references
\addtocontents{toc}{\protect\setcounter{tocdepth}{2}}
\bibliographystyle{plain}
\bibliography{references.bib}
\end{document}